\documentstyle{amsppt}
\input epsf
\NoBlackBoxes 
\pagewidth{6.5 in}
\mathsurround=1pt
\define\today{\the\month.\the\day.\the\year}
\topmatter
\title
	Combinatorically Prescribed Packings\\
	And Applications to Conformal and Quasiconformal Maps
\endtitle
\rightheadtext{Combinatorically Prescribed Packings}
\author Oded Schramm \endauthor
\affil Princeton University \endaffil
\abstract{
The Andreev-Thurston Circle Packing Theorem
is generalized to packings of convex bodies in
planar simply connected domains.
This turns out to be a useful tool
for constructing conformal and quasiconformal mappings
with interesting geometric properties.
We attempt to illustrate this with a few results about
uniformizations of finitely connected planar domains.
For example, the following variation of a theorem
by Courant, Manel and Shiffman is proved and generalized.
If $G$ is an $n+1$-connected bounded planar domain,
$H$ is a simply connected bounded planar domain,
and $P_1,P_2,\dots,P_n$ are (compact) planar
convex bodies, then sets $P_j'$ can be found so that
$G$ is conformally equivalent to
$H-\cup_{j=1}^n P_j'$, and each $P_j'$ is either a point,
or is positively homothetic to $P_j$.
}\endabstract
\endtopmatter
\define\fffffig#1:#2:#3:#4:{%
\midinsert
\centerline{\epsfysize=#4\epsfbox{#3.ps}}
{\noindent \centerline{{\bf Figure #1. #2}}}\endinsert}
%
\redefine\leq{\leqslant}
\redefine\geq{\geqslant}
\define\ri{\therosteritem}
\countdef\secnum=1
\secnum=0
\countdef\fignum=2
\fignum=0
\define\newsec{\global\advance\secnum by 1 \fignum=0}
\define\fig{\the\secnum.\the\fignum}
\define\newfig{\global\advance\fignum by 1 \fig}

\redefine\t#1{{{#1}^*}}
\define\br{\overline}
\define\dedg#1#2{\langle{#1},{#2}\rangle}
\define\tri#1#2#3{\langle{#1},{#2},{#3}\rangle}
\define\G#1#2{LS_{#1}^{#2}}
\define\g#1#2{RS_{#1}^{#2}}
\define\m{{\goth M}}
\define\f{{\goth F}}
\define\gp{\gamma'}
\define\gam{\gamma}
\define\Gam{\Gamma}

\define\eps{\epsilon}

\document
\baselineskip=20pt 

\newsec\heading
	1. Introduction 
\endheading

\subheading{Packing Theorems}
Say we are given a finite packing of circles in the plane or in the sphere.
The nerve of the packing is a planar graph
(without multiple edges or loops)
which describes the combinatorics of the packing.
The vertices of the nerve correspond to the circles in the packing,
and the edges correspond to pairs of circles which touch.
A reinterpretation by Thurston\footnote
{Given at the International Symposium in Celebration of the Proof
of the Bieberbach Conjecture. Purdue University, March 1985.
See also \cite{Th}.}
of a theorem  of Andreev (\cite{An1}, \cite{An2}) states that for any planar graph
there exists a packing of circles
whose nerve is the given graph.
We will refer to this fact as the Circle Packing Theorem.

We generalize the Circle Packing Theorem by proving the
existence of packings involving more general figures.
For example, we prove the following.
Suppose that a planar graph is given,
together with a correspondence which assigns to each
vertex of the graph a (compact) smooth planar convex body.
Then it is possible to modify these convex bodies by homotheties,
so that they form a packing whose nerve is the given graph. 
The Circle Packing Theorem is obtained
when one takes all the convex bodies to be disks.
The proof of this generalization does not rely on previously known
proofs of the Circle Packing Theorem, and is self contained,
except for the use of Brower's Fixed Point Theorem.

It is possible to use the methods of this paper to prove a generalization
regarding packing of balls of Riemannian metrics:
Given a planar graph together with a correspondence which assigns to each
vertex of the graph a Riemannian metric on the sphere, it is always possible
to find a packing on the sphere with the graph as its nerve so that each packed
set is a ball in the Riemannian metric corresponding to the vertex associated to
the set.  See \cite{Sch1}, \cite{Sch2}.

Since the completion of this work, the author has developed two additional
approaches which yield similar generalizations of the Circle Packing Theorem
\cite{Sch2}, \cite{Sch3}.  The advantages of these newer methods are
that they give shorter proofs, and that one obtains sharp uniqueness
statements.  The main benefit of the technique presented here is that
it can produce packing theorems for which no reasonable uniqueness
statement holds, and is in that sense more general.

\subheading{Conformal Mappings}
Thurston conjectured that the Circle Packing Theorem
is related to the Riemann Mapping Theorem in that
circle packings can be used to yield approximately
conformal mappings.
This conjecture was later proved by Rodin and Sullivan in the
charming paper \cite{R-S}.
Basically, the idea is that conformal
mappings are characterized by the property that their differentials
take circles to circles. 
To be more specific, consider a simply connected bounded planar domain $G$.
Look at a hexagonal packing in the plane with small circles; that is,
an infinite packing of small circles, all having the same size,
so that each one touches six others.  From this infinite packing select
those circles that are contained in $G$, and call the resulting packing $H$.
From the Circle Packing Theorem it follows that
one can find a circle packing $H'$ contained in the unit disk, together with
a correspondence which assigns to every circle $c$ in $H$ a circle
$c'$ in $H'$, in such a way that touching circles in $H$ correspond
to touching circles in $H'$, and that a boundary circle in $H$ (that is,
a circle which doesn't touch six other circles of $H$) corresponds to 
a circle of $H'$ which touches the unit circle (from the inside).  
The correspondence
$c\to c'$ tends to the conformal map from $G$ onto the unit disk,
as the size of the circles in the hexagonal packing tends to zero,
provided appropriate normalizations are made.

Pursuing this idea, we use our ability to pack convex
shapes and the techniques of \cite{R-S} to prove and generalize
the following theorem about conformal mappings of multiply connected domains.

\proclaim{9.1 Theorem}
Let $G$ be an $n+1$-connected bounded domain in $\Bbb C$
which is obtained from a simply connected region $H$
after $n$ disjoint compact
connected sets $F_1,F_2,\dots,F_{n}$ have been removed from it.
Assume that none of the $F_j$ is a point. (To avoid trivialities).
Then for every planar bounded simply connected region $H'$,
and for every list 
$P_1,P_2,\dots P_{n}$ of $n$ compact convex sets which are not points,
there are disjoint sets
$P_1',P_2',\dots,P'_{n}$ contained in $H'$, with
each $P_j'$ (positively) homothetic to $P_j$,
such that $G$ is conformally equivalent to 
$
H'-\cup_{j=1}^{n}P_j'
$
and the boundary of each $F_j$ corresponds to
the boundary of the respective $P_j'$ under such a conformal equivalence.

In the case where the boundaries of $H$ and $H'$ are simple
closed curves, $\gamma $ and $\gamma'$ respectively, if three distinct
counterclockwise ordered points are chosen $z_1,z_2,z_3\in \gamma $,
and similarly  $z_1',z_2',z_3'\in \gamma '$,
then such $P_j'$-s and such an equivalence can be chosen
to satisfy the additional requirement that
each $z_j$ corresponds to $z'_j$ when the conformal
equivalence is extended continuously to $\gamma ,\gamma '$.
\endproclaim

The first part of this theorem follows from the work of
Courant, Manel and Shiffman \cite{C-M-S}.
Their results are obtained by variational techniques, and
the proof given here is very different.  (See the introduction
in \cite{Si} for references and a short survey of related results).
In  Theorem~9.1, instead
of prescribing the sets $P_j'$ up to homothety, one can use more
general prescriptions.  For example, if one is given a $C^1$
foliation of
a neighborhood of $\overline{H'}$ by curves, then the
theorem still holds if $P_1'$ is required to be some subarc of
one of the leaves in the foliation.  Meanwhile, the other $P_j'$ can
be prescribed as in the theorem, or by the same foliation, or by other
such foliations.
Though this general statement follows from the techniques of this paper
(packings of Riemannian balls), we prove it here only for the case where
the leaves of the foliation are straight line segments.
The proof for foliations by curves will appear in a subsequent paper.

\subheading{Quasiconformal Maps}
Circle packings have also been studied in connection to quasiconformal maps.
Z.~X.~He \cite{He2} used the Circle Packing Theorem to construct
solutions to the Beltrami equation.  However, for the study of quasiconformal
maps it is more natural to use packings of ellipses than packings of circles.
The reason for this is that while conformal maps 
are characterized by the fact that their differential takes circles
to circles, quasiconformal maps are (almost) characterized by the
fact that their differential takes an ellipse of bounded eccentricity
to an ellipse of bounded eccentricity.
In the same manner that circle packings can be used to construct
conformal maps, packings of ellipses should yield quasiconformal
maps $f$ that their differential at a point $z$ takes a specified
ellipse to another ellipse specified up to homothety.
(The shapes of the ellipse in the domain and of the ellipse
in the range are assumed to depend continuously on $z$ and $f(z)$.)
This idea, for which I am indebted to Peter Doyle, can be
made to work in many situations.  One example for this is a generalization
of the above Theorem~9.1 to quasiconformal maps satisfying a generalized
Beltrami equation, which we state and prove in Section~9.

\subheading{Organization of the paper}
In Section 2 we give some simple definitions and introduce notations.
Section 3 is a sketch of the technique used here to prove the
packing theorems.  In Sections 4, 5 some preparatory work is done,
and more definitions are given.  At the end of Section 5, a general
packing theorem is stated, the Monster Packing Theorem.  The subsequent
packing theorems are directly derived from the Monster Packing Theorem. 
Sections 6 and 7 are devoted to the 
proof of this theorem.  Section 8 contains the results concerning  packings
of convex planar sets, while in Section 9 these results are applied to
yield conformal and quasiconformal mappings of multiply connected domains.

\subheading{Acknowledgements}
I am deeply grateful to my advisor Bill Thurston, to Peter Doyle,
and to Richard Schwartz, for teaching me much about circle-packings,
for many discussions, and for giving 
much valuable advice concerning this work.  This paper is an edited version
of my Ph.~D.\ thesis \cite{Sch1}.

\newsec\heading
	2. Basic Definitions and Notations
\endheading

For us, \lq the plane\rq\ will be the complex plane, $\Bbb C$, and
\lq the sphere\rq\ will be the Riemann sphere,
$\hat{\Bbb C}={\Bbb C}\cup\{\infty\}$.

Let $Q_1,Q_2$ be two closed sets in the sphere.
We will say that $Q_1$ and $Q_2$ {\it touch}, if they intersect, but 
the interior of $Q_1$
is disjoint from $Q_2$, and the interior of $Q_2$ is disjoint from $Q_1$.
If they intersect, but do not touch, we shall say that they {\it intersect
nontrivially}.
A {\it packing\/} is a finite indexed collection of closed, connected, nonempty 
sets, such that there are no nontrivial intersections.
A packing $Q=(Q_i:i\in I)$ in the sphere.
is {\it degenerated}, 
if there is a point in common to three of the sets $Q_i$,
or if one of the $Q_i$-s contains only a single point.
Generally, when we mention a packing, it is assumed to be nondegenerated;
otherwise we will say specifically that it may be degenerated.

The {\it nerve\/} of the packing $Q=(Q_i:i\in I)$ is the (abstract)
graph whose vertex set is $I$ and which has a (single) edge
$\dedg ij$ joining two
distinct vertices $i,j$ if and only if the sets $Q_i,Q_j$ touch.
If each of the packed sets $Q_i$ is topologically a disk, then the
nerve is a planar graph; that is, it can be embedded in the plane. 
To see this, choose an arbitrary interior point,
$v_i$, in each of the sets $Q_i$, and choose a simple curve $\gam_{i,j}$ in
$Q_i\cup Q_j$ to join $v_i$ and $v_j$ whenever $Q_i$ and $Q_j$ touch.
It is clearly possible to do this in such a way that the various
edges $\gam_{i,j}$ are disjoint, except at the vertices, and
this defines a planar embedding of the nerve.
(Here the fact that the packing is nondegenerated is used).



Let $T$ be a triangulation of the sphere (without loops or multiple edges),
and let $G$ be its 1-skeleton as an abstract graph.  Then $G$
is a planar graph, and it is not hard to see that $T$ can be reconstructed
from $G$ (up to homeomorphism).  This is proved by noting that 
a triangle in $G$ is the boundary of a triangle of $T$ if and only
if it does not separate $G$.
For this reason, we will be a little sloppy, and will not distinguish
between a triangulation and its 1-skeleton.  For example, we may say that
the nerve of a packing is a triangulation.

When trying to prove the existence of a packing with a specified nerve,
and satisfying some conditions on the shapes of the sets,
it is generally sufficient to consider nerves which are 
triangulations, because if one is given a planar graph which is not a
triangulation, then additional vertices may be added,
and edges may be adjoined to these new vertices, connecting
them to each other and to preexisting vertices, to yield a triangulation.
Thus, we will be dealing here mostly with triangulations.
When we use the term {\it triangulation\/} it is implicitly assumed
that we mean a triangulation of the sphere without loops and without
multiple edges.  Furthermore, the triangulations are assumed to
be oriented, and so terms like \lq clockwise\rq\ can be used.

If $Q(x)$ is a compact set in the sphere, dependent on a parameter $x$,
then when we refer to
the limit of $Q(x)$ as $x\to y$, we will mean the limit in the
sense of the Housdorff metric;
likewise, continuity of $Q(x)$ is to be understood
with respect to this metric.

A {\it convex body\/} is a bounded compact convex set (in
the plane) which has interior points.  A {\it homothety\/}
is a transformation of the form $z\to az+b$, where $a,b$ are
constants, $a>0$ is real, and $b\in\Bbb C$.  Two sets are
{\it homothetic}, if there is a homothety taking one set onto
the other.

\newsec\heading
	3. An Outline of the Packing Technique
\endheading
In this section we try to sketch the technique used to prove the
existence of packings on which our results are based.
We do this by restricting ourselves to a simplified situation,
and giving a more or less complete proof for that case.
The core of the proof consists of an $n$-dimensional argument,
and there are also many, mostly simple, two dimensional arguments.
Hopefully this section will give the reader an idea of the core
of the proof, without entangling him or her in the myriad of
details. 

Instead of examining a general triangulation, we shall now consider
only triangulations which have a spiral:

\definition{3.1 Definition}
Let $T$ be an (oriented) triangulation of the sphere, and
suppose that $v_1,v_2,\dots,v_s$ is a Hamiltonian path in $T$.
(That means that $v_1,v_2,\dots,v_s$ are distinct, that they are all
the vertices of $T$,
and that $v_j$ neighbors with $v_{j-1}$, for $j=2,3,\dots,s$).
Further suppose that the edge $\dedg{v_j}{v_k}$ lies on the left
side of the path $v_{j-1}\to v_j\to v_{j+1}$,
whenever $k>j+1,\ j>1,$ and $v_k$ neighbors with $v_j$; and that
$\dedg{v_j}{v_k}$ lies on the right side of $v_{j-1}\to v_j\to v_{j+1}$,
whenever $k<j-1,\ j<s,$ and $v_k$ neighbors with $v_j$.
Then $v_1,v_2,\dots,v_s$ is called a {\it counterclockwise spiral\/} in $T$.
See Fig.~\newfig.
\enddefinition
\fffffig \fig:A triangulation with a counterclockwise spiral.:\fig:2.0in: 

Say we are given a triangulation of the sphere $T$, together with
a counterclockwise spiral $v_1,v_2,v_3,\dots,v_s$ in $T$.
Let $\gam$ be a smooth simple closed curve in the plane $\Bbb C$,
and let $Q_1\subset\hat{\Bbb C}={\Bbb C}\cup\{\infty\}$ be the closure in
$\hat{\Bbb C}$
of the region determined by $\gam$ which contains $\infty$.
Suppose that $Q_2,Q_3$ are two circles which touch each other,
and each one of them touches $Q_1$.
Our objective in this section is to establish the following proposition.

\proclaim{3.2 Proposition}
Let $Q_1,Q_2,Q_3,$ and $T$ be as above.
There are circles $Q_4,Q_5,\dots,Q_s$ which, together with
$Q_1,Q_2,Q_3$, form a packing whose nerve is isomorphic to $T$.
The exhibited isomorphism will take the vertex of the nerve which
corresponds to $Q_j$ to the vertex $v_j$ of $T$, for $1\leq j\leq s$.
\endproclaim

Set $C_j{\overset{\text{def}}\to=}Q_j,\ j=1,2,3$.
We examine the configurations of circles
$C_4,C_5,\dots,C_s$ satisfying, for $j=4,5,\dots,s$, the following conditions.
(Here we allow some of these circles to degenerate to points).
\roster
\item"{(3.1)}" $C_j$ touches the circle $C_{j-1}$.
\item"{(3.2)}" $C_j$ also touches $\cup_{t=1}^{j-2}C_t$.
(This also means that the interiors are disjoint).
\endroster

We parametrize this set of configurations in the following way.
Take the $(s-3)$-dimensional cube
$$
K=\{(x_4,x_5,\dots,x_s):0\leq x_j\leq 1,\ j=4,5,\dots,s\}.
$$
Each point $x\in K$ will determine such a configuration of circles
$C_4(x),C_5(x),\dots,C_s(x)$. 
Set $C_j(x){\overset{\text{def}}\to=}C_j,\ j=1,2,3,\ x\in K$.
The definition
of $C_j(x)$ will be inductive; that is, $C_j(x)$ depends on $C_k(x),\ k<j$.
Let $4\leq j\leq s$, and suppose that the circles
$C_k(x),\ 4\leq k<j$ have been placed already.  The coordinate $x_j$ will
help us place $C_j(x)$.  Look at $C_{j-1}(x)$.  There is a unique point,
say $p_{j-1}=p_{j-1}(x)$, in which $C_{j-1}(x)$ touches $C_{j-2}(x)$.
Let $p_j=p_j(x)$ be the point
whose angular distance from $p_{j-1}$ is $2\pi x_j$,
measured counterclockwise along the boundary of $C_{j-1}(x)$.
The circle $C_j(x)$ is chosen as the largest circle which touches $C_{j-1}(x)$
at $p_j$ and its interior is disjoint from $\cup_{k=1}^{j-2}C_k$.
(If no such circle exists, set $C_j(x)=\{p_j(x)\}$.
In the special case where $C_{j-1}(x)=\{p_{j-1}(x)\}$, we take
$p_j(x)=p_{j-1}(x),\ C_j(x)=C_{j-1}(x)$.)
It is not too difficult to see that this procedure defines
the circles $C_j(x)$, that this definition
satisfys \ri{3.1}, \ri{3.2}, and that the circles $C_j(x)$ depend
continuously on $x$.

We wish to demonstrate the existence of a point $x\in K$ such that
the configuration corresponding to it is a packing with the required
nerve.
Take some $j,\ 4\leq j\leq s$.  If $x_j=0$, or $x_j=1$, then all the circles
$C_k(x),\ k\geq j$ degenerate and are the point of contact of
$C_{j-1}(x)$ with $C_{j-2}(x)$.  If $x_j$ is close to zero,
then all the circles $C_k(x),\ k\geq j,$ are trapped in the
funnel between $C_{j-1}(x)$ and $C_{j-2}(x)$. (See Fig.~\newfig).
If $x_j$ is close to one, they are trapped in the other
funnel between $C_{j-1}(x)$ and $C_{j-2}(x)$.
The first funnel is located near the \lq right\rq\ side of $C_{j-2}(x)$,
and the second is located near the \lq left\rq\ side.
Here, and in the following, the {\it left side\/} of the circle
$C_i(x),\ (2\leq i<s),$ is that part of the boundary
of $C_i(x)$ which extends counterclockwise from the point $p_{i+1}(x)$,
where $C_{i+1}(x)$ touches $C_i(x)$, to $p_i(x)$, the point
where $C_i(x)$ touches $C_{i-1}(x)$, including $p_i(x)$ and $p_{i+1}(x)$.
\fffffig \fig:Some circles in the funnel.:\fig:4in:

When trying to find an $x$ which gives our elusive packing,
we see that if $x_j$ is close to zero, it is too small,
because all the circles $C_t(x),\ t\geq j,$ are in one of the funnels,
and if it is close to one, it is too large,
because all the circles $C_t(x),\ t\geq j,$ are in the other funnel.
One hopes that between zero and one there is some value
which is right for $x_j$.
Of course, this must be done simultaneously for $j=4,5,\dots,s$.
This hints that some kind of $(s-3)$-dimensional
Intermediate Value Theorem must be used, and this is exactly what we will do.

Denote the 0-faces and 1-faces of $K$ by
$$
F^0_j=\{x\in K: x_j=0\},\quad F^1_j=\{x\in K: x_j=1\},\qquad j=4,5,\dots,s.
$$
For each $j=4,5,\dots,s$, let $h_j$ be the least $h\geq j$ such that
$v_{j-2}$ has an edge in $T$ with $v_h$.
(To see that $v_{j-2}$ neighbors with some vertex $v_t,\ t\geq j$,
look at the edge $\dedg {v_{j-2}}{v_{j-1}}$ of $T$.
Since $T$ is a triangulation, there is some vertex, say $v_t$,
which completes
this edge to a triangular 2-cell of $T$
lying to the left of $\dedg {v_{j-2}}{v_{j-1}}$.
Now, $t>j-1$, since $v_1,v_2,\dots,v_s$ is a spiral.
So we see that $v_{j-2}$ neighbors with some vertex $v_t,\ t\geq j$,
and this shows that the definition above is OK).
Define
$$
\align
&K_j=\{x\in K:\cup_{k=j}^{h_j}C_k(x)\text{ touches the
left side of }C_{j-2}(x)\}-F^0_j,\\
&L=K-\cup_{j=4}^sK_j.
\endalign
$$
The proof of Proposition~3.2 consists of showing that
$\br L\cap(\cap_{j=4}^s \br K_j) \neq \emptyset$,
and that for a point $y$ in this intersection, the sets
$C_1(y),\dots,C_s(y)$ form a packing as required.

To see that the intersection $\br L\cap(\cap_{j=4}^s \br K_j)$
is not empty, the following facts will be used.
\roster
\item"{(3.3)}" $K=\cup_{j=4}^s K_j\cup L$.
\item"{(3.4)}" $K_j\supset F_j^1$, for $j=4,5,\dots,s$.
\item"{(3.5)}" The relative interior of $\cup_{i=4}^{j-1}K_i$ contains
     $\cap_{i=4}^s\br K_i\cap F_j^0$, for $j=4,5,\dots,s$.
    (\lq The relative interior\rq\ means the interior relative
    to $K$).
\endroster
Statement \ri{3.3} is obvious.
Statement \ri{3.4} holds,
because if $x_j=1$, then $C_j(x)$ contains $p_{j-1}$, the intersection point
of $C_{j-1}(x)$ and $C_{j-2}(x)$.
(Actually, in this case $C_j(x)=\{p_{j-1}\}$).
Statement \ri{3.5} will be proved shortly.

Let us consider a point $y$ in $\cap_{i=4}^s \br{K_i}$, and
let $j$ be in the range $4\leq j\leq s$. 
As $y\in\br{K_j}$,
there is a sequence of points $x\in K_j$ converging to $y$.
Since the definition of $K_j$ implies that for these $x$,
$\cup_{k=j}^{h_j}C_k(x)$ touches the left side of $C_{j-2}(x)$,
we know that 
$$
\cup_{k=j}^{h_j}C_k(y)\text{ touches the left side of }C_{j-2}(y),\tag 3.6
$$
by continuity.

What happens if any of the circles $C_i(y)$ degenerate to a point?
This certainly cannot happen when $i<4$.
Assume that this does happen for $i\geq 4$,
and that $i$ is the minimal index such that $C_i(y)$ is a point.
We shall see now that $y$ is in the relative interior of 
$K_i$.
Since $C_i(y)$ is a point, suppose that $C_i(y)=\{p\}$.
Then, by the construction, $C_k(y)=\{p\}$ for $k\geq i$.
Taking $j=i$, the union in \thetag{3.6} contains only the point $p$,
and therefore $p$ is on the left side of $C_{i-2}(y)$.
Since $C_{i-1}(y)$ and $C_{i-2}(y)$ are touching circles, we conclude
that $p$ is their point of tangency.  This means $y_i=0$ or
$y_i=1$.

Examine the case $y_i=1$ first.  The sets $C_1(y),\dots,C_{i-3}(y)$
are smooth; they are not points.
None of these sets contains the point $p$, because
$C_{i-1}(y)$ and $C_{i-2}(y)$ do, and it is impossible for three
smooth sets to touch at a point.
So there is a positive distance from $p$ to $\cup_{t=1}^{i-3}C_t(y)$.
This implies that, for $x$ close enough to $y$, the circle $C_i(x)$
touches the left side of $C_{i-2}(x)$, because we know that it touches
$\cup_{t=1}^{i-2}C_t(x)$, and $x_i$ is close to $1$. 
But this means that $y$ is in the relative interior of $K_i$.

For the case $y_i=0$, we use a similar argument.  For $x\notin F_i^0$
close to $y$
all the circles $C_k(x),\ k\geq i$ are in the funnel between $C_{i-1}(x)$
and the right side of $C_{i-2}(x)$.  They therefore cannot reach the
left side of $C_{i-2}(x)$, and these $x$ are not in $K_i$.  This
contradicts $y\in \br {K_i}$, and shows that
$y_i\neq 0$.

Now we are in a position to prove statement~\ri{3.5}.
Suppose that $x$ is a point in $\cap_{i=4}^s\br K_i\cap F_j^0$.
As $x_j=0$, it is easy to see that $C_j(x)$ consists of a single
point.  Therefore, we may substitute $x$ for $y$ in the above analysis,
to conclude that $x$ is in the relative interior of $K_i$, for some $i<j$.
This verifys \ri{3.5}.

In Section~6 we will use Brower's Fixed Point Theorem
to see that the facts \ri{3.3}, \ri{3.4}, \ri{3.5} imply that the
intersection $\br L\cap(\cap_{j=4}^s \br K_j)$ is not empty. 
Other topological arguments can be used instead.  

Having convinced ourselves that 
$\br L\cap(\cap_{j=4}^s \br K_j) \neq \emptyset$,
we pick some point $y$ in this intersection.
It remains to verify that the sets $C_1(y),\dots,C_s(y)$
form a packing as required.

Since we are now assuming that $y\in \br L$, in addition to 
$y\in \cap_{i=4}^s \br{K_i}$, the discussion above shows that none
of the sets $C_j(y)$ degenerates to a single point, because otherwise
we would have $y\in\text{relative interior }K_i$, for some $i$,
contradicting $y\in \br L$.

The following lemma gathers combinatorial properties of the spiral
$v_1,v_2,\dots,v_s$ that we will need below.

\proclaim{3.3 Lemma}
Let $j$ be an index in the range $3<j\leq s$, and
let $k_j$ be the smallest index $k$ such that
$v_j$ neighbors with $v_k$ in $T$. (See Fig.~\newfig).
Then
\roster
\item $k_j<j-2$,
\item $v_{j-1}$ neighbors with $v_{k_j}$,
\item if $k_j<m\leq j-2$, then $v_m$ has no neighbor in the range
$v_t,\ m+1<t<j$.
\endroster
\endproclaim
\fffffig \fig:A section of a spiral.:\fig:2.8in: 

\demo{Proof}
Since $T$ is a triangulation, there is some vertex, say $v_k$,
which forms a triangle together with the edge
$\dedg {v_{j-1}}{v_j}$, and this triangle, $\tri {v_k}{v_{j-1}}{v_j}$,
lies to the right of the edge $\dedg {v_{j-1}}{v_j}$.

$k\neq j-2$, because otherwise the edge $\dedg{v_{j-2}}{v_j}$
would emerge from the right side of the path
$v_{j-3}\to v_{j-2}\to v_{j-1}$.
The edge $\dedg {v_k}{v_{j-1}}$ approaches the path
$v_{j-2}\to v_{j-1}\to v_{j}$ from the right, and therefore $k<j-1$.
Using $k\neq j-2$, we get $k<j-2$.
This gives \ri1, since $k_j\leq k$.

To prove \ri2, we will show that $k=k_j$.
If $k=1$, this is obvious.  Assume $k>1$.
The edge $\dedg {v_k}{v_{j-1}}$ is to the right of the
path
$v_{j-2}\to v_{j-1}\to v_{j}$,
and to the left of the path
$v_{k-1}\to v_{k}\to v_{k+1}$.
Therefore the closed path 
$$
v_{k}\to v_{k+1}\to \dots\to v_{j-2}\to v_{j-1}\to v_k
$$
separates $v_j$ from $v_{k-1}$, and therefore from all the
vertices $v_t,\ t\leq k-1$. 
So $v_j$ does not neighbor with these vertices, and $k=k_j$.
Thus \ri2 is established.

The proof of \ri3 is similar, and will not be given here.
\qed\enddemo

We have used the fact that $y\in\br {K_j}$ to prove \thetag{3.6}.
What can we deduce from $y\in\br L$?
Note that $y\notin\cup_{i=4}^sF_i^0$, because we have shown that
the circles $C_i(y)$ do not degenerate.
Since $y\in\br L$, there is a sequence of $x$ in $L-\cup_{i=4}^sF_i^0$
converging to $y$.
Take such an $x$.  Because $x\notin \cup_{i=4}^jK_i$, the circle $C_j(x)$
cannot touch the left side of a circle $C_{i-2}(x)$,
when $4\leq i\leq j\leq h_i$.
Let $k_j$ be the smallest index $k$ such
that $v_j$ neighbors with $v_k$ in $T$.
If $k_j<i-2\leq j-2$,
then $4\leq i\leq j$ holds, and
$j\leq h_i$ follows from \ri3 of Lemma~3.3.
So, restating our previous conclusion,
the circle $C_j(x)$ does not touch the left side of any of the circles
$C_t(x),\ k_j<t\leq j-2$.
But, by the construction, $C_j(x)$ touches $\cup_{t=1}^{j-2}C_t(x)$.
We therefore see that $C_j(x)$ either touches 
$\cup_{t=1}^{k_j}C_t(x)$, or touches the right side of one of the
circles $C_t(x),\ k_j<t\leq j-2$.
Using continuity, we see that the same holds for $y$ instead of $x$; that is,
for $j=4,5,\dots,s$,
$$
C_j(y)\text{ touches }
\left(
(\cup _{t=1}^{k_j} C_t(y))
\cup
(\cup_{t=k_j+1}^{j-2} \text{the right side of } C_t(y))
\right)
.\tag 3.7
$$

\thetag{3.6} and \thetag{3.7} will be the only tools that we will
use to show that the touchings dictated by the triangulation $T$
actually occur in our packing
$C_1(y),C_2(y),\dots,C_s(y)$.
This is done in the following lemma.

\proclaim{3.4 Lemma}
If $v_i$ neighbors with $v_q$ in $T$ and $i<q$,
then $C_q(y)$ touches the left side of $C_i(y)$. 
\endproclaim
\demo{Proof}
The lemma clearly holds when $q=i+1$.
The proof for $q>i+1$ will be by induction on $q$.
The case $q\leq 3$ being obvious, we take some $q,\ 3<q\leq s$,
and assume that the lemma holds for smaller values of $q$.

{\smc Case 1}, $i=k_q$ (notation of Lemma 3.3):
We know that $v_{q-1}$ neighbors with $v_i$, from \ri2 of Lemma 3.3.
The induction hypotheses then implies that $C_{q-1}(y)$
touches the left side of $C_i(y)$.
Let $p$ be a point in this intersection.
Consider the closed curve which travels counterclockwise along the boundary
of $C_i(y)$ from $p$ to the first point of intersection
with  $C_{i+1}(y)$, then travels counterclockwise along the
boundary of $C_{i+1}(y)$, till it reaches $C_{i+2}(y)$, and so on,
until, while walking along the boundary of $C_{q-1}(y)$, it encounters
$p$ again.
This curve determines a region $R\subset \Bbb C-\cup_{t=1}^{q-1}C_t(y)$,
whose boundary consists of part of the left side of $C_i(y)$, the left sides of
each of the circles $C_t(y),\ i<t<{q-1}$, 
and part of the boundary of $C_{q-1}(y)$.
See Fig.~\newfig.
\fffffig \fig:The region $R$.:\fig:4in: 

Taking $j=q$ in \thetag{3.6}, shows that one of
the circles $C_t(y),\ t\geq q,$ touches the left
side of $C_{q-2}(y)$, and therefore (because $i=k_q<q-2$) is in $\br R$.
But if one of these is in $\br R$, then all of them are, because they are
connected.  In particular, $C_q(y)\subset\br R$.
This shows that $C_q(y)$ cannot touch the sets $C_t(y),\ t<k_q$,
and cannot touch the right side of the sets $C_t(y),\ k_q\leq t < q-1$.
By \thetag{3.7} with $j=q$, $C_q(y)$ must therefore touch the
left side of $C_{k_q}(y)=C_i(y)$, as required.

{\smc Case 2}, $k_q<i<q-1$:
In this case $q=h_{i+2}$, by \ri3 of Lemma 3.3.
Looking at \thetag{3.6} with $j=i+2$, we see that it is sufficient to
show that none of the circles $C_t(y),\ i+2\leq t<q$, touches the left
side of $C_i(y)$.  
 Assume that this is not the case: assume that $i+2\leq m<q$, and that
$C_m(y)$ touches the left side of $C_i(y)$.  A contradiction will
be achieved using an argument very similar to the one appearing in 
case 1.

Consider the region $R'$ of $\Bbb C-\cup_{t=1}^mC_t(y)$ whose boundary
consists of part of the left side of $C_i(y)$, the left sides of
each of the circles $C_t(y),\ i<t<m$, and part of the boundary of $C_m(y)$.
Taking $j=m+1$ in \thetag{3.6} shows that one of
the circles $C_t(y),\ t\geq m+1$, is in $\br{R'}$.
But if one of these is in $\br{R'}$, then all of them are, 
in particular $C_q(y)$.  This contradicts
what we have established in case 2, that $C_q(y)$ touches 
$C_{k_q}(y)$.
This contradiction completes the discussion of case 2 and
the proof of the lemma.
\qed\enddemo

Clearly, Proposition 3.2 follows from Lemma 3.4, and so its proof is
complete.

Reviewing what we have done in this section will reveal
that we have used only the very general properties
of circles.  In fact, there is no difficulty to generalize the discussion to
wider classes of sets.  The harder part is to get rid of the assumption
that the triangulation considered has a spiral, but the difficulty
is mostly technical.

\newsec\heading
	4. The Tree
\endheading

Let $T$ be an arbitrary, but fixed,
triangulation of the sphere, 
and let $\dedg ab$ be an edge in $T$.
We shall denote the vertices of $T$ by $V$.
It is the purpose of this section to construct a special spanning
tree in $T$, called a {\it right oriented depth first search tree \/},
and to derive those properties which we will need later.
This spanning tree will take the role of the spiral from
the previous section.

To construct the tree we imagine the vertices as belonging to a dynasty,
and the tree will describe their descendence hierarchy.
The founder of the dynasty is $a$.  He produces the
child $b$ and dies. 
There are several
rules that govern the lives of the vertices and the evolution of the
dynasty.
Every vertex lives only once.
Except $a$, every vertex has one parent; $a$ has none. 
Just before he dies, a vertex may have children which are 
neighbors of it in $T$. 
The vertices live a sad and lonely life. 
The loneliness rule says that in every
given moment, and in every connected component of 
$T-\{\text{dead vertices}\},$ there is exactly one live vertex.

Let's follow the life of a typical vertex, $v$.  As he is
born and looks around him, what does he see?  He sees 
a circle of vertices consisting of his
neighbors in $T,$ some of them dead, and some not yet born. 
All the unborn neighbors must eventually become descendents
of $v,$ because $v$ is the only live vertex in its connected
component of nondead vertices.
On the other hand, if $u$ is a dead neighbor of $v$, then,
by the same reasoning as above with $v$ replaced by $u$,
we see that $u$ is an ancestor of $v$.

So the circle of neighbors of $v$ is partitioned into arcs containing
ancestors of it and arcs containing descendents.
The (maximal) arcs consisting of descendents will 
be called \lq\lq the arcs\rq\rq.
It is obvious that two descendents included
in the same arc are in the same connected component of
the unborn vertices in $T$.
The converse is also true: if $u,w$ are two neighbors of $v$
that are in the same connected component of the unborn vertices,
then they are actually in the same arc.
To see this, look at a simple (not self-intersecting),
path of unborn vertices joining $u$
and $w$.  When  $v$ is joined to this path one obtains
a simple closed path.  This closed path 
separates $T$ into two components, because
the topological space underlying $T$ is a sphere.
The component that doesn't
contain the root $a$ must contain only unborn
vertices.  The neighbors of $v$ that are in this component
connect $u$ and $w,$ thus showing that $u$ and $w$ are
in the same arc.

In order to preserve the loneliness rule
when $v$ dies, it is necessary and sufficient that
$v$ will have exactly one child in every
connected component of the unborn vertices that is contained 
in his connected component of $T-\{\text{dead vertices}\}$.
By the discussion above, this means that $v$ has one
child in each arc of unborn vertices surrounding him.
The last rule governing the evolution of the dynasty
is that every vertex $v\neq a,$
has the most clockwise
vertex in every such arc as a child.  This is meaningful,
since $a$ is the only vertex having the
unborn vertices completely surrounding him (Every
other vertex has a parent). 

It is immediate that
this procedure defines a spanning tree, $R,$ of $T,$
with root $a,$ the edges of $R$ being 
$\{\dedg fs :f \text{ is the parent of }s\}$.  
$R$ defines a partial order, $\leq,$
on the vertices, in the following manner.  We say that
$u\leq v,$ if $v$ is a descendent of $u$ or $v=u$.
Equivalently, $u\leq v$ if and only if $u$ is on the
shortest tree-path from $v$ to the root $a$.
Since the partial order $\leq$ is induced, in this
manner, by a tree, it has the property that for each
$i\in V$ the set $\{j\in V:j\leq i\}$ is linearly
(that is, totally) ordered.

We have nearly completed the proof of the following proposition.
 
\proclaim{4.1 Proposition}
There is a spanning tree $R$, in $T$, with root $a$,
so that $R$ and
the order, $\leq$, induced by it satisfy:
\roster
\item  $b$ is the only child of $a$.
\item  For each edge in the triangulation, its vertices, $e,f$,
are comparable.
Either $e\leq f,$ or $f\leq e$.
\item  If $v\neq a$ is a vertex in $T$ and $A$ is a maximal arc
consisting of descendents of $v$ neighboring with $v,$ then all the vertices
in $A$ are also descendents of the most clockwise vertex in $A$, which
is a child of $v$.
\item  Let $u,w$ be descendents of $v$ neighboring with $v$,
then their meet (their latest common ancestor), is $v$, or
neighbors with $v$.
\endroster
\endproclaim

When $v\in T-\{a\}$ we will denote by $v-1$ the parent of $v$
in the tree.  Similarly $v-2$ denotes the grandparent of $v,$
in case $v\neq a,b$.  Do not be mislead by the notation ---
$v-1=u-1$ does not imply $u=v$.
The tree constructed above will be called a {\it right oriented
dfs tree for\/} $T$ {\it rooted at\/} $\dedg ab$.  
Dfs stands for depth-first-search.

\demo{Proof}
Only \ri 4 remains to be justified.
Let $s$ be the meet of $u,w$ and assume $s\neq v$.
Let $u'$ be the child of $s$ which is $\leq u$,
and similarly for $w'$ and $w$.  We have $u'\neq w'$.
Let $A'$ be the arc of neighbors of $s$, extending
between $u'$ and $w'$, which does not contain $s-1$.
From \ri 3, we know that there is some ancestor $t$
of $s$ in $A'$.

Look at the cycle $C$ composed of the (shortest) tree-path from
$s$ to $u$, the edges $\dedg vu, \dedg vw$ and the
tree-path from $w$ to $s$.  
This cycle separates $T$ into two components, because
the topological space underlying $T$ is a sphere.
($C$ may separate the graph of $T$ into more than two
components, but we take into account the two-cells as well).
Looking in the neighborhood of $s$, we 
notice that this cycle separates $t$ from
$s-1$, if they are not in $C$.

Now look at the circle of neighbors of $v$.  As $u,w$ descend from
the same child of $v$, one of the arcs extending between them
must consist exclusively of descendents of $v$, by \ri 3.
Denote this arc by $A$.  Again from \ri 3, we see that in $A$
$v$ has no children.
From this we conclude that the component of $T-C$ which corresponds
to $A$ contains only vertices which are $\geq s$, because all the vertices
in $C$, except $v$, are $\geq s$.  
Therefore $s-1$ cannot be in this component, and neither can $t$.
Since they cannot be together in the other component, at least one 
of them is in $C$.  Both are $<s$, so the only possibility is
that one of them is $v$.  In any case $v$ neighbors with $s$.
\qed
\enddemo

Statement \ri 3 of the Proposition has the following useful consequence.
If $u,\ v$ are neighboring
vertices in $T$ and $u-1>v>a$, then the edge $\dedg vu$ grows out
of the left side of the (directed) tree path from $a$ to $u-1$.
(Compare to Definition 3.1.)
Such an edge will be called a \lq left edge\rq.
So we see that all the edges in $T$ are either tree-edges or left edges.

There is a unique vertex in $T$, which we will denote by $c$,
so that $\tri a b c$ is a clockwise triangular 2-cell of $T$.
From \ri3 of the proposition it follows that $b=c-1$, and that 
all the vertices of $T$, other than $a,b,c$ are descendents of $c$;
that is, $v>c$ for every $v\in V-\{a,b,c\}$.

The following lemma provides more detail about the structure of $(T,\leq)$.

\proclaim{4.2 Cycle Lemma}
Let $i>c$ be a vertex in our triangulation $T$.  
Let $I=\{j\in V:j\geq i\}$ and let $H$ be the collection of
vertices in $T$
bounding $I$.
(The vertices in $H$ are those vertices of $T-I$ that 
have neighbors in $I$). 
Then $H$ forms a cycle in $T$. (Fig.~\newfig).
Furthermore, we have,
\roster
\item
The vertices in $H$ are linearly ordered (by our usual order induced
from the tree). 
Denoting the vertices of $H$ by 
$v_0<v_1<v_2<\dots <v_e$, we have $v_e=i-1$.
\item 
The edges $\dedg {v_0}{v_e}$
and $\dedg {v_{j-1}}{v_j},\ j=1,2,\dots e$
occur in $T$.
If $0<j<e$ and $v_{j+1}<v<i$, then the edge $\dedg {v_j}v$
does not occur in $T$.
\item 
$\tri {i-1}i{v_0}$ forms a clockwise triangle in $T$.
For $0\leq j\leq e-1$, there is a vertex $i'\in I$, so that
$\tri {v_j}{i'}{v_{j+1}}$ forms a clockwise triangle in $T$.
\item
The vertices of $I$ lie on the left side of the cycle
$v_0\to v_1\to\dots v_e\to v_0$.
\endroster
\endproclaim

{
  \midinsert
  \epsfysize=3.5in
  \centerline{\epsfbox{\fig.ps}}
  \centerline{{\bf Figure \fig. The cycle $H$.}}
  \narrower\narrower\narrower\narrower
  \narrower\narrower\narrower\narrower
  \vskip 0in
  \noindent
  \it An arrow, $u\to w$ has been placed between
  two vertices $u, w$, whenever $u=w-1$.
  \vskip 0in
  \endinsert
}

\demo{Proof} 
Let $u$ be any vertex in $H$.  It has a neighbor in $I$, say $j$.
We can walk along the boundary of $H$ and $I$ in the following manner.

Walking forward:  
$u\notin I$ and therefore $u-1\notin I$, so $u$ has some neighbor
which is not in $I$.  
Looking at the circle of
neighbors around $u$,  we move clockwise, starting at $j$,
until a vertex not in $I$ is reached.  Denote this vertex by $u'$,
and the vertex in $I$ just before it denote by $j'$.  This is one
step.  Now the next step can be taken if we put $u'$ in 
place of $u$ and $j'$ in place of $j$.
Note that $\tri u{j'}{u'}$ forms a clockwise triangle in
$T$.

Let us check the case where $u>u'$.
In that case
\ri 3 of Proposition 4.1 implies that $j'$ is a child of $u$,
($j'>u$ is true in any case).
Since $j'\in I,\ u\notin I$, we conclude that
$j'=i,\ u=i-1$.

What happens if we start walking from the pair $i-1,i$?  
Call the vertex of $H$ that is reached in the first step
$v_0$, and the vertex of $I$ reached call $i_0$.
Let the vertices of $H$ that are reached in subsequent
steps be $v_1,v_2,\dots$, and the vertices of $I$
that are reached in subsequent steps be $i_1,i_2,\dots$.
Denote by $e$ the first index $\geq 0$ so that
$v_e > v_{e+1}$;
i.e., $v_0<v_1<\dots<v_e,\ v_e>v_{e+1}$.
From the discussion above we may conclude that
$i_{e+1}=i, v_e=i-1, v_{e+1}=v_0$.  

It is clear that the first part of \ri 2 is satisfied.
Also \ri 3 has been established, with $i'=i_{j+1}$.
To verify \ri 4, note that the cycle in \ri 4 separates $T$
into a component lying to the left of it and a component lying to the
right of it.  As $I$ is connected, it must be contained in one of these.
From \ri 3, we see that this must in fact
be the component to the left.

To see that there are no other vertices in $H$, we will point
out that walking forward is invertible --- one can walk backwards.
Walking backward is the same as forward, 
we only replace \lq clockwise\rq\
with \lq counterclockwise\rq.
Walking forward and then backward from the pair $u,j$ returns us 
to $u$.  $j$ may perhaps be replaced by $h$, another vertex 
in $I$, but this will have no consequences.  \lq Have
no consequences\rq\ here means that walking forward or
backward from $u,j$ gives the same result as from
$u,h$.

Starting from any adjacent pair $u_0\in H,\ i_0\in I$, we
define $u_j,i_j$ by walking forward
in the obvious manner.  There will be some
index $f$ so that $u_f>u_{f+1}$.  Then
we have $u_f=i-1,\ i_{f+1}=i$.  Therefore
all the vertices in $H$ are
those obtained by walking backward from $i-1,i$.  This
gives $H=\{v_0,v_1,\dots v_e\}$.

It remains to prove the second part of assertion \ri 2.
Assume to the contrary that $0<j<e$, $v_{j+1}<v<i$ and the edge
$\dedg{v_j}v$ occurs in $T$.
Look at the cycle $C$ formed from the tree-path from $v_j$
to $v$ by adjoining the edge $\dedg{v_j}v$.
Let this cycle be oriented in the direction of the tree-path from $v_j$
to $v$, and let the region bounded by the left side of this
cycle be denoted by $R$.
There is some vertex in $I$, say $k$, which neighbors with $v_{j+1}$.
Since the edge $\dedg {v_{j+1}}k$ is a left edge,
we see that $k$ is necessarily in $R$.
This shows that $i$ is also in $R$, because the tree path from $i$
to $k$ cannot cross $C$.  Because $\dedg {v_j}v$ is a left edge,
$v_j-1$ is not in $R$, and this shows that $v_0$ is not in
$R$, contradicting the fact that $i$ neighbors with $v_0$.
The contradiction establishes the second part of assertion \ri 2,
and completes the proof of the lemma.
\qed
\enddemo

\newsec\heading
	5. Monsters
\endheading

In Section 3 we have defined a mapping which assigns to each
point in a $\operatorname{Card}(V)-3$ cube a configuration of
circles in the plane.  That mapping is a particular case of 
what we call a {\it monster}.
 
\definition{5.1 Definition}
Let $T$ be a triangulation of a 2-sphere without multiple edges
or loops, and let $\tri abc$ be a clockwise triangle in $T$.
Let $\leq$ be the partial order on $V$, the vertices of $T$, induced by a
right-oriented dfs tree rooted at the edge $\dedg ab$.  (As we have seen,
$a<b<c$ are the three smallest vertices.)  Set $J=V-\{a,b,c\}$.
A {\it monster\/} for $(T,\leq)$ is a mapping $\m$, which assigns to
each point $x$ in the 
$\operatorname{Card}(J)$-dimensional cube 
$K=[0,1]^J$ a configuration of sets $(\m_i(x): i\in V)$ in
the 2-sphere, $\hat{\Bbb C}$, such that the following conditions
are satisfied.
\roster
\item"{(M1)}" $\m_i(x)$ is topologically a closed disk or a point, for
	$i\in J,\ x\in K$.  When $i\in \{a,b,c\}$ and $x\in K$, 
	$\m_i(x)$ is topologically a closed disk.
\item"{(M2)}" The mapping is continuous.  This means that for every
	$i\in V$ the set
	$\m_i(x)$ tends to $\m_i(y)$, as $x$ tends to $y$ in $K$. 
\item"{(M3)}" The interior of $\m_i(x)$ is disjoint from every
	$\m_j(x),\ j<i$
	(Note that when $i,j$ are not comparable, overlaps are permitted).
\item"{(M4)}" $\m_i(x)$ touches $\m_{i-1}(x)$, if $i\neq a$.
	Furthermore, there are points $p_i=p_i(x)$ in the intersection
	$\m_i(x)\cap \m_{i-1}(x)$, which depend continuously on $x$.
	The point $p_i(x)$ will be called the {\it base\/} of $\m_i(x)$.
\item"{(M5)}" $\m_i(x)$ touches one of the sets $\m_j(x),\ j<i-1$, if
	$i\neq a,b$.
\item"{(M6)}" For each $i\in J$ and $x\in K$,
	$p_i(x)=p_{i-1}(x)$ if and only if either $x_i=0$, $x_i=1$ or
	$\m_{i-1}(x)$ is a point.
	Furthermore, if $x$ varies in $K$ and $x_i$ tends to zero,
	then the diameter of the counterclockwise arc from $p_{i-1}(x)$
	to $p_i(x)$ along the boundary of $\m_{i-1}(x)$ tends to
	zero, and if $x_i$ tends to $1$, then the diameter of the
	clockwise arc from $p_{i-1}(x)$ to $p_i(x)$ along the boundary of
	$\m_{i-1}(x)$ tends to zero. 
%
\item"{(M7)}" No-three-touch rule:  When $i<j<k,\ x\in K$ and 
	$\m_i(x)\cap \m_j(x)\cap \m_k(x)\neq \emptyset$, 
	then one of these three sets is a point.
\endroster
\enddefinition
Condition \ri{M6} is worthy of some explanation.  It is meant to 
describe the motion of $\m_i(x)$ around $\m_{i-1}(x)$.  It says that
$\m_i(x)$ circles counterclockwise around $\m_{i-1}(x)$,
relative to $p_{i-1}$,
as $x_i$ increases from $0$ to $1$.  The reason for the slight
technicality comes from the fact that the set $\m_{i-1}(x)$ may
vary during this process.

Note that \ri{M7}, the no-three-touch rule,
is automatically satisfied if the sets $\m_i(x)$,
when they don't degenerate to points,
have smooth boundaries, or don't have \lq\lq sharp\rq\rq\ corners.
This is because of \therosteritem{M3}.
Sometimes we will write $\m_i$ in place of $\m_i(x)$; likewise $p_i=p_i(x)$.

One can just see the terrible monster swinging its arms in sheer
rage, the tentacles causing a frightful hiss, as they rub against
each other.  Luckily, the following theorem comes to the rescue.

\proclaim{5.2 Monster Packing Theorem}
Let $\m$ be a monster for $(T,\leq)$ and let $K,V$ be as above.  
Then there is a point
$y\in K$, such that the sets $(\m_i(y):i\in V)$ form a packing 
with nerve $T$.
\endproclaim

A few lemmas will prepare for the proof of this theorem,
which will be completed in Section 7.
In Section 8 monsters consisting (mostly) of convex sets will
be constructed, and existence results for packings of convex
sets will follow from the Monster Packing Theorem.

\newsec\heading
	6. The Decomposition
\endheading
Let us fix a monster $\m$ and all the relevant data:
$T,\ \leq,\ V,\ J=V-\{a,b,c\},\ K$ etc.
In this section we decompose the cube
$K=[0,1]^J$ according to the touchings occurring between the
$\m_i(x),\ i\in V,\ x\in K$.  We will see that the
closures of the subsets forming this decomposition intersect.
In the next section it will be proved that a point in this
intersection corresponds to a packing as required for the Monster
Packing Theorem.

Denote the 0-faces and 1-faces of $K$ by
$$ F^0_i=\{x\in K: x_i=0\},\quad F^1_i=\{x\in K: x_i=1\}.$$

Let $j<i$ be vertices in $V-\{a\}$.  
Suppose, for instance, that the sets $\m_i$ and $\m_j$
touch, for a given $x\in K$.  To say only that they touch
does not describe the situation well.  They can touch
in two distinct ways (disregarding degenerate situations
in which at least one of these sets happens to be a
point).
$\m_i$ can touch either the right side or the left
side of $\m_j$.  
Here, as usual, \lq right\rq\ and
\lq left\rq\ are relative terms --- they depend
not only on $j$ but also on $i$.  To make the distinction
we will need the following notation.

\definition{6.1 Notation}
Let $a\neq j<i\in V,\ x\in K$.
Look at the \lq snake\rq\ of sets
$\cup_{j\leq k <i}\m_k(x)$.  The {\it left side\/} of the boundary
of this snake will be denoted by $\G ji$, or sometimes by
$\G ji(x)$.  More precisely,
we define
$\G ji=\cup_{j<k\leq i}\G{k-1}k$, where $\G {k-1}k$, the
left side of $\m_{k-1}$, is the 
part of the boundary of
$\m_{k-1}(x)$ extending clockwise from $p_{k-1}$ to $p_k$, 
including the endpoints $p_{k-1}$ and $p_k$.
(If $\m_{k-1}(x)$ happens to degenerate to a point $p$, then
we set $\G {k-1}k=\{p\}$).
When $j=a$ we define $\G ji$ similarly, but include the whole
boundary of $\m_a$.

The {\it right side\/}, $\g ji$, is defined similarly, with
\lq clockwise\rq\ replaced by \lq counterclockwise\rq.
\enddefinition

Note that the sets $\G ji$ and $\g ji$ vary continuously
as $x$ varies in $K$.  Note also that in general
$\G ji\neq \G j{i'}$ when $i,i'$ are distinct children
of $j$.

For every $i\in J$ we define a subset $K_i\subset K,$ in the
following manner.  Look at $i$ from the direction of
$i-1 $.  $i$ is the most clockwise vertex in an arc $A$ of
descendents of $i-1 $ neighboring with $i-1 $. 
The neighbor of $i-1$ bounding $A$ in the counterclockwise
direction is an ancestor of $i-1$.  It will be called the
{\it godparent\/} of $i$, and will be denoted $g_i$.
See Fig.~\newfig.
There is a $v\geq i$
that is a neighbor of $g_i$ in $T,$
namely the most counterclockwise vertex in $A$.
From \ri4 of Proposition 4.1, we see that of all vertices $v\geq i$
there is a unique minimal one $h_i$ which is a neighbor of $g_i$
in $T$. (Possibly $h_i=i$).
$h_i$ is called the {\it heir\/} of $i$.
Now set
$$
\align
&K_i=\{x\in K:
	(\cup_{i\leq j\leq h_i} \m_j)
    \cap
	\G {g_i}{i-1}\neq\emptyset\}-F_i^0,\\
&L=K-\cup_{i\in J} K_i.
\endalign
$$
Our decomposition of $K$ is
$$
K=\cup_{i\in J} K_i\cup L.
$$

\fffffig \fig::\fig:2.0in: 

\proclaim{6.2 Intersection Lemma}
$$\cap_{i\in J} \br {K_i}\cap \br L\neq \emptyset.$$
\endproclaim

The proof of this lemma will be given below, following
some preparation.
A point $y\in \cap_{i\in J} \br{K_i}\cap \br L$
will henceforth be called an {\it intersection point}.
Our method of proving the Monster Packing Theorem will be to
show that any intersection point $y$ gives a packing,
$(\m_i(y),\ i\in V),$ 
whose nerve is $T$.

We do not know much about the structure of our decomposition of
$K$.
However, we will be able to obtain some information about
the boundary behavior of the decomposition.
The Intersection Lemma will follow from this information
and a topological argument (Brower's Fixed Point Theorem).

\proclaim{6.3 Boundary Behavior Lemma}
\roster
\item $K=\cup_{j\in J} K_j\cup L$.
\item $K_j\supset F_j^1$ for $j\in J$.
\item For every $j\in J$, $F_j^0\cap (\cap_{i\in J}\br{K_i})$
  is contained in the relative interior of $\cup_{i<j}K_i$.
  (\lq The relative interior\rq\ means the interior relative to $K$).
\endroster
\endproclaim
\demo{Proof}
\ri1 follows from the definition of $L$.
For \ri2 recall that $\m_j=\{p_{j-1} \}$ when $x_j=1$, 
and $p_{j-1}$ is in the left side of $\m_{j-2}$.

\ri3 will require some work.  Let $j\in J$, and suppose that
$x\in F_j^0\cap (\cap_{i\in J}\br K_i)$.
Since $x_j=0$, $p_j(x)=p_{j-1}(x)$ follows from \ri{M6}.
Now, if $\m_{j-1}(x)$ and $\m_{j-2}(x)$ do not degenerate to 
points, then $\m_j(x)$ does, because $p_{j-1}\in\m_{j-1}\cap\m_{j-2}$,
and because of the no-three-touch rule, \ri{M7}.
So we see that one of the sets $\m_i(x),\ i\leq j$ degenerates to a point.

We will need to study the situation where some of the sets degenerate,
anyway,  
so let us consider a more general situation.

\proclaim{6.4 Claim}
Let $y$ be a point in
$\cap_{i\in J}\br K_i$, and suppose that one of the sets $\m_i(y),\ i\in V,$
degenerates to a point.
Let $i$ be a minimal vertex, 
in the order induced by the tree, so that
$\m_i(y)$ is a point.
Then $i\neq a,b,c$, $y\notin F_i^0$, and
$y$ is in the relative interior of $K_i$.
\endproclaim
\demo{Proof}
By \ri{M1} $i\neq a,b,c$.
Suppose that $\m_i(y)=\{p\}$. 
From \ri{M4} and \ri{M5}, we know that $p$ is a point where
$\m_{i-1}(y)$ intersects with $\m_k(y)$, for some $k<i-1$.
Because of the minimality of $i$, the sets $\m_{i-1}(y)$
and $\m_k(y)$ are not points.
If $t$ is a child of $i$ ($i=t-1$), then $\m_t(y)$ contains
$p$ and therefore necessarily $\m_t(y)=\{p\}$, by the no-three-touch
rule \ri{M7}.  Continuing, we see that $\m_k(y)=\{p\}$, for all $k\geq i$.

Since $y\in \br K_i,$ arbitrarily close to $y$ there are
points $y'\in K_i$.  For these $y'$,
$$
	(\cup_{i\leq k\leq h_i} \m_k)
    \cap
	\G{g_i}{i-1}\neq\emptyset.\tag 6.1
$$
We will examine three cases.

{\smc Case 1}, $y\in F^0_i$:
We wish to see that this cannot happen.
In this case we have $p=p_{i-1}(y)$, by \ri{M6}.
Look at a point $y'\in K_i$ very close to $y$.
Provided
$y'$ is close enough to $y$, the
$\m_k(y'),\ k\geq i,$ will still be very small and very close to
$p_{i-1}(y')$, by continuity. 

Look at the union of the bodies corresponding to descendents
of $i$, including $i$:
$$
D(x)=\cup_{k\geq i} \m_k(x).
$$
By \thetag{6.1}, $D(y')$ intersects
$\G{g_i}{i-1}(y')$.
Let us first see that $D(y')\cap\G{i-2}{i-1}(y')=\emptyset$.

Since $y'\in K_i,$ $y'_i\neq 0,$ and since $y'$ is close to $y$,
$y'_i$ is small. 
Because of \ri{M6}, this means
that the base of $\m_i(y')$, $p_i(y')$, is
a little counterclockwise from $p_{i-1}(y')$ on the boundary
of $\m_{i-1}(y')$. 
So $\m_i(y')$ and $\G{i-2}{i-1}(y')$ 
are on different sides of $\m_{i-2}(y') \cup \m_{i-1}(y')$.
(See Fig.~\newfig).
In other words, for $D(y')$ to touch $\G{i-2}{i-1}(y')$ it must either
go around $\m_{i-1}(y')$ or go around $\m_{i-2}(y')$, or sneak
between them. (The interiors of the $\m_k,\ k\geq i$
are disjoint from $\cup_{k<i} \m_k$).
But for $y'$ sufficiently close to $y$, all the
bodies composing $D(y')$ are very small in relation to
$\m_{i-1}(y')$ and $\m_{i-2}(y'),$ by continuity.  Therefore,
they cannot manage to go around.
Because of the no-three-touch rule, it is also
obvious that the chain $D(y')$ can't sneak between
$\m_{i-1}(y')$ and $\m_{i-2}(y')$.
\fffffig \fig::\fig:3in: 

So we conclude that $D(y')$ doesn't intersect $\G{i-2}{i-1}(y')$.
But \thetag{6.1} holds for $y',$ so $D(y')$ intersects
$\G{g_i}{i-2}(y')$.
However, this also can't happen when $y'$ is close enough
to $y,$ because then $D(y')$ is contained in an arbitrarily
small neighborhood of $p=p_{i-1}(y)$, and $p_{i-1}(y)$ is not in
$\cup_{g_i\leq k<i-2}\m_k(y)$,
by the no-three-touch rule, since $p_{i-1}\in \m_{i-1}\cap \m_{i-2}$.
This contradiction shows that $y\notin F^0_i$.

{\smc Case 2}, $y\notin F^0_i\cup F^1_i$:
In this case we have $p\neq p_{i-1}(y)$,
by \ri{M6}.
We again consider $y'$ close to $y$ for which 
\thetag{6.1} is satisfied.
Since intersection of compact sets
is a closed condition, and $\m_k(z),\ \G tk(z)$ 
depend continuously on $z\in K$, \thetag{6.1} also 
holds for $y$.  
Because $D(y)=\m_i(y)=\{p\},$ this leads us to conclude
that $p$ is a point where $\m_{i-1}(y)$ intersects with
$\G{g_i}{i-1}(y)$.
There is a small neighborhood of $p$ which intersects only
one of the sets $\m_k(y),\ k<i-1$, and intersects only the
left side of its boundary.  This is true because of the
no-three-touch rule and because $p\in \m_{i-1}(y),\ p\neq p_{i-1}(y)$.
(The place where the left and right sides of the boundary meet
is a place where two consecutive sets intersect).
For any $z,$ $\m_i(z)$ intersects with some $\m_k(z),\ 
k<i-1$.  
So for $z\in K$ close enough to $y,$ $\m_i(z)$ still intersects
$\G{i-1}{g_i}(z)$.
We infer that \thetag{6.1} holds for all $z\in K$ sufficiently
close to $y$.  Therefore $y$ is in the interior of
$K_i$ relative to $K$, as we wanted to show.

{\smc Case 3}, $y\in F^1_i$:  This case is similar.
When $y'$ is close to $y$, $\m_i(y')$ is either $\{p_{i-1}\}$ or is
small and located very close to $p_{i-1}$ and a little clockwise from 
it along the boundary of $\m_{i-1}(y')$, by \ri{M6}.  
But $\m_i(y')$ must touch some $\m_k(y'),\ k<i-1$,
and we see, by the same arguments used above,
that this $k$ can only be $i-2$ (for $y'$ close
enough to $y$), and, in fact, $\m_i(y')$ must touch the left 
side of $\m_{i-2}(y')$.  This implies $y'\in K_i$ and $y$ is
in the interior of $K_i$ relative to $K$,
as needed.
\qed\enddemo

\demo{Proof of the Boundary Behavior Lemma, Conclusion}
We were in the midst of proving \ri3.  We may take
$y=x$ in Claim~6.4, and choose $i$ to be a 
minimal vertex $\leq j$ so that $\m_i(y)$ is a single
point.  The Claim then says that $i$ is actually $<j$,
and that $x$ is in the relative interior of $K_i$.
Since $x$ is an arbitrary point in $F_j^0\cap (\cap_{i\in J}\br K_i)$,
\ri3 holds, and the proof is complete.
\qed\enddemo

\proclaim{6.5 Topological Lemma}
Let $K_j',\ j\in J,$ and $L'$ be subsets of $K$,
and assume that
$$
\gather 
  K=\cup_{j\in J} K'_j\cup L',\tag 6.2\\
  {\align   K_j'&\supset F_j^1,\qquad j\in J,\tag 6.3\\
          L'&\supset \cup_{j\in J}F_j^0.\tag 6.4  \endalign}
\endgather
$$
Then
$$
\cap_{j\in J}\br{K_j'}\cap\br {L'}\neq \emptyset.
$$
\endproclaim

This lemma is well-known.
However, since it can easily be deduced from 
Brower's Fixed Point Theorem, 
we will give a proof here.

\demo{Proof}
In order to apply 
Brower's Fixed Point Theorem we define a continuous
function $f:K\to K$.  For $j\in J$ let $e^j$ be the
vector whose $j$-th coordinate is 1 and whose other
coordinates are 0.  Let $e^0$ be the vector 
with all coordinates $-1$.
Define $d_j(x)$ to be the
distance of $x$ from $K'_j,$\ $j\in J,$ and $d_0(x)$
to be the distance of $x$ from $L'$.  Set
$$f(x)=x+\sum _{j\in J\cup\{0\}} d _j(x)e^j.$$
Let $x$ be a point in $K$.  Since $K'_j\supset F_j^1,$
$d_j(x)\leq 1-x_j,$ so the $j$-th coordinate of
$f(x)$ is not greater than 1.  This, together with a 
similar argument applied to $L'$ and the zero-faces
$F^0_j,$ shows that $f$ maps $K$ into itself.
Using Brower's Fixed Point Theorem, 
let $y$ be a point in $K$ satisfying
$f(y)=y$.

Because $K=\cup_{j\in J} K'_j\cup L'$, at least one
of the $d_j(y)$ is zero. 
But $f(y)=y$ implies
$$0=\sum _{j\in J\cup\{0\}} d _j(y)e^j.$$
This shows that all the $d_j(y)$ are zero,
because every proper subset of 
$\{e^j:j\in J \text{ or }j=0\}$
is linearly independent.
Thus $y$ is in $\br{ L'}$ and in the $\br {K'_j},\ j\in J,$
as required.
\qed
\enddemo

\demo{Proof of the Intersection Lemma}
The results of the Boundary Behavior Lemma are almost
the hypotheses of the Topological Lemma.  Only small
modifications are necessary.

Define
$$
\align
&U_j=(\text{the relative interior of }\cup_{i<j}K_i)-F_j^1,\qquad j\in J,\\
&K_j'=K_j-U_j,\qquad j\in J,\\
&L'=L\cup(\cup_{i\in J}F_i^0).
\endalign
$$
\thetag{6.3} holds, because of \ri2 of the Boundary Behavior Lemma 6.3,
and because $U_j$ is disjoint from $F_j^1$.
It is clear that \thetag{6.4} holds.
To prove \thetag{6.2} consider a point $x\in K$.  If $x \notin\cup_{j\in J}
U_j$, then certainly $x\in \cup_{i\in J} K'_i\cup L'$, because
of \ri1 of the Boundary Behavior Lemma.  Otherwise, if $j$ is a minimal
element of $J$ such that $x\in U_j$, then $x\in\cup_{i<j}K_i$,
and also $x\in\cup_{i<j}K_i'$.
Thus \thetag{6.2} holds.

Using the Topological Lemma 6.5, we see that
$\cap_{i\in J}\br{K_i'}\cap\br {L'}\neq \emptyset$.
Let $y$ be a point in this intersection.
Obviously, $y\in \br{K_j},\ j\in J$, since $K_j\supset K_j'$.
Therefore, \ri3 of the Boundary Behavior Lemma and $y\notin U_j,\ j\in J$,
imply that $y\notin F_0^j,\ j\in J$.
So we see that $y$ is also in $\br L$, and the proof is complete.
\qed
\enddemo
\enddemo

\newsec\heading
	7. Intersection Point(s)
\endheading

We shall now fix an intersection point
$y\in\br L \cap(\cap_{i\in J}\br K_i)$.
The Intersection Lemma 6.2 guarantees the existence of such a $y$.
Note that, since $y\in\br L$, and $L$ is disjoint
from $\cup_{i\in J}K_i$, Claim 6.4 implies that $y$ is not
in any of the zero-faces, $F_i^0$, and none of the sets $\m_j(y)$  
degenerates to a point.  We will show that these sets actually form
the required packing.

\proclaim{7.1 Contacts Lemma}
Let $i\in J$, let $v_0<v_1<\dots<v_e$ be as in the Cycle Lemma~4.2,
and recall the definitions of $g_i,h_i$.
Then
$$
\cup_{i\leq j\leq h_i}\m_j(y)\text{ touches }\G{g_i}{i-1}(y),
$$
and
$$
\m_i(y)\text{ touches } \cup_{j<v_1}\m_j(y) \cup \g {v_1}{i-1}(y).
$$
\endproclaim

\demo{Proof}
Take some $i\in J$.  Arbitrarily close to $y$ there are
points $y'\in K_i$.  For these points $y'$, 
$$
	(\cup_{i\leq j\leq h_i} \m_j)
    \cap
	\G{g_i}{i-1}\neq\emptyset\tag{7.1}
$$
holds.
Since the sets in \thetag{7.1} are compact, and vary continuously as
$y'$ varies in $K$, \thetag{7.1} also holds for $y$.
This verifies the first assertion.

Because $y\in\br L$, from Claim 6.4 it follows that $y$ is
not in any of the zero faces, therefore $y$ is in the closure of
$L-\cup_{j\in J}F_j^0$.
Let $z$ be a point in $L-\cup_{j\in J}F_j^0$.
Like in the proof of our first assertion, it is sufficient to show
that
$\m_i(z) \cap \left(\cup_{j<{v_1}}\m_j(z) \cup \g {v_1}{i-1}(z)\right)
\neq\emptyset$.
From \ri{M5}, we know that $\m_i(z)$ touches $\cup_{j<i-1}\m_j(z)$,
and therefore it suffices to show that
$\m_i(z)$ doesn't touch $\G {v_1}{i-1}(z)$.
To see this we will use the fact $z\notin \cup_{j\in J}K_j$.

Let $t$ be some index in the range $0<t<e$, and let $j$ be the child
of $v_{t+1}$ which is $\leq i$.  ($v_{t+1}=j-1,\ j\leq i$).
From \ri3 of the Cycle Lemma 4.2 and \ri3 of Proposition 4.1, it follows
that $g_j=v_t$.  From the second part of \ri2 of the Cycle Lemma,
it then follows that $i\leq h_j$.  $z\notin K_j$ now implies that
$\m_i(z)$ does not touch $\G {g_j}{j-1}=\G {v_{t+1}}{v_t}$.
Since $\G {v_1}{i-1}=\cup_{0 < t<e}\G{v_t}{v_{t+1}}$,
we know that $\m_i(z)$ cannot touch $\G {v_1}{i-1}(z)$,
and this completes the proof.
\qed
\enddemo

From now on we will write $\m_i$ instead of $\m_i(y)$, $(i\in V)$, since
we will no longer use the sets $\m_i(x),\ x\neq y$.
Similarly, we will use $\G ij$ in place of $\G ij(y)$, and
$\g ij$ in place of $\g ij(y)$.

\proclaim{7.2 Inductive Lemma}
Let $i\in V-\{a\}$.  For each $v<i-1$ neighboring with $i$ in $T$,
$\m_i$ touches the left side of $\m_v$.
\endproclaim

Before proving the lemma, We shall introduce the following notation.
Suppose that $u$ and $v$ are vertices of $T$ satisfying $u-1>v$ 
and $\m_u$ touches the left side of $\m_v$.
Then there is a region in the complement of $\cup_{v\leq j\leq u}\m_j$
whose boundary consists of part of $\partial \m_u$ together with $\G vu$ 
excluding a part of $\partial\m_v$.  See Fig.~\newfig.
We will denote the closure of this region by $R_v^u$. 
(This region may be non-unique, if $\m_u\cap\m_v$ contains more
than a single point.  In that case, we choose $R_v^u$ to 
denote one of those regions.)
If $j<v$ or $j>u$, and $\m_j$ intersects $R_v^u$, then
$\m_j$ is contained in $R_v^u$ This follows from the
no-three-touch rule, because $\m_j$ cannot intersect
nontrivially with $\cup_{v\leq j\leq u}\m_j$.  (Again see Fig.~\fig). 
If $v\neq a$, then $R_v^u$ does not contain the base of $\m_v$, 
which is in $\m_{v-1}$.  Therefore $R_v^u$ is disjoint from $\m_{v-1}$.
Continuing, we see that $R_v^u$ is disjoint from all the sets
$\m_j,\ j<v$.  Similarly, if $j>u$ and $\m_j\subset R_v^u$,
then all the sets $\m_k$ with $k$ comparable to $j$ and
$k>u$ are contained in $R_v^u$.
We will use these facts in the following.
\fffffig \fig:The region $R_v^u$.:\fig:3in:  

\demo{Proof}
As the name suggests, the lemma will be proved by induction.
It is obviously true for the case $i=a,b,c$.  
Take some $i\in J$.
We will assume the lemma holds when $i$ is replaced by any
vertex $<i$, and prove that it holds for $i$.

Let $v_0<v_2<\dots<v_e=i-1$ be as in the Cycle Lemma 4.2.
By the inductive hypotheses, $\m_{v_e}$ touches the left side of
$\m_{v_0}$.
From the first part of the Contacts Lemma, we know that
$\cup_{i\leq j\leq h_i}\m_j$ touches $\G{v_{e-1}}{i-1}$.
This shows that one of the sets $\m_j,\ j>i$ intersects $R_{v_0}^{v_e}$
and therefore is contained in $R_{v_0}^{v_e}$. 
But if one of the sets $\m_j,\ j>i$ is contained in 
$R_{v_0}^{v_e}$, then all of them are, for they are connected.
In particular, $\m_i\subset {R_{v_0}^{v_e}}$.

Our next objective is to show that $\m_i$ touches the left side of $\m_{v_0}$.
Since $\m_i$ is in ${R_{v_0}^{v_e}}$,
it cannot possibly touch any set $\m_j,\ j<v_0$,
and it cannot touch $\g{v_0}{v_e}$.
The second part of the Contacts Lemma
says that 
$\m_i$ touches $\cup_{j<v_1}\m_j\cup \g {v_1}{i-1}$,
and therefore, to show that 
$\m_i$ touches the left side of $\m_{v_0}$,
it is sufficient to prove that $\m_i$ does not touch the left
side of any of the sets $\m_k,\ v_0<k<v_1$. 
This is trivial if $v_1-1=v_0$, so assume $v_1-1>v_0$.
Let $s$ be the child of $v_1$ which is $<i$. 
See Fig.~\newfig.
Since $i$ is not in the region of $T$ which is bounded by the
left side of the tree-path from $v_0$ to $v_1$ and the edge $\dedg {v_0}{v_1}$,
neither is $s$.
Let $v$ be
the ancestor of $v_1$ which is one vertex clockwise from $s$ among
the neighbors of $v_1$. 
$v$ cannot be $>v_0$, because it is impossible for the edge $\dedg vs$ to
go out of the right side of the tree-path from $v_0$ to $v_1$;
and $v\neq v_0$, because $v = v_0$ would mean that there are two distinct
edges between $v_0$ and $v_1$.  
So we conclude that $v<v_0$.
By the induction hypotheses, $\m_s$ touches $\m_v$. This shows that
$\m_s$ is not in $ {R_{v_0}^{v_1}}$, and therefore $\m_i$ isn't
there either.
So in any case $\m_i$ does not touch the left side of any of the sets
$\m_k,\ v_0<k<v_1$. We conclude that $\m_i$ touches the left
side of $\m_{v_0}$.
\fffffig \fig:A portion of $T$.:\fig:2.5in:  

Any vertex $v<i-1$ which neighbors with $i$ must be one of the
vertices $v_0,v_1,\dots,v_{e-1}$.  Assume that $v_t$ neighbors with
$i$ ($0\leq t\leq e-1$).  We must show that $\m_i$ touches the left
side of $\m_{v_t}$.  We have done this already for $t=0$, so assume
that $t>0$.
Let $j$ be the child of $v_{t+1}$ which is $\leq i$.
($v_{t+1}=j-1,\ j\leq i$).  In the proof of the Contacts Lemma
we have observed that, in this setting, $g_j=v_t$ and $i\leq h_j$.
But $i\geq h_j$, because $i$ neighbors with $v_t$.
So $h_j=i$.
From the Contacts Lemma we know that
$\cup_{j\leq k\leq h_j}\m_k$ touches $\G{g_j}{j-1}=\G{v_t}{v_{t+1}}$.
Therefore, to verify that $\m_i$ touches the left side of $\m_{v_t}$,
it is sufficient to show that none of the sets $\m_k,\ j\leq k<i$ touches
$\G{v_t}{v_{t+1}}$, and that $\m_i$ cannot touch the left
side of any of the sets $\m_k,\ v_t<k<v_{t+1}$.

We will show now that $\m_i$ cannot touch the left
side of any of the sets $\m_k,\ v_t<k<v_{t+1}$.
If $v_{t+1}-1=v_t$ there is nothing to show.  Otherwise, by the inductive
hypotheses, $\m_{v_{t+1}}$ touches the left side of $\m_{v_t}$.
Since we have shown that $\m_i$ touches $\m_{v_0}$, it follows that
$\m_i$ is not in ${R_{v_t}^{v_{t+1}}}$, and therefore
$\m_i$ does not touch the left side of any of the sets
$\m_k,\ v_t<k<v_{t+1}$, as was required.

We shall show now that none of the sets $\m_k,\ j\leq k<i$ touches
$\G{v_t}{v_{t+1}}$.
Since $\m_i$ is not in ${R_{v_t}^{v_{t+1}}}$, the same is true
for all the sets $\m_u,\ u\geq j$, and therefore they too
cannot touch the left side of any of the sets $\m_k,\ v_t<k<v_{t+1}$.
It remains to show that none of the
sets $\m_u,\ j\leq u<i$ touches the left side of $\m_{v_t}$.
Assume to the contrary that $j\leq u<i$ and $\m_u$ touches the left side of
$\m_{v_t}$.
Let $r$ be the child of $u$ which is $\leq i$.  It is easy to see that
$g_r>v_t$.  $\m_i$  is not in ${R_{v_t}^u}$, because it touches
$\m_{v_0}$.  This implies that none of the sets $\m_k,\ k\geq r$ are
there, in contradiction to the Contact Lemma which says that
$\cup_{r\leq k\leq h_r}\m_k$ touches $\G{g_r}{r-1}$.
This contradiction completes the proof that $\m_i$ touches
the left side of $\m_{v_t}$, and the proof of the lemma.
\qed\enddemo

\demo{Proof of the Monster Packing Theorem}
We will first show that $(\m_i:i\in V)$ is a packing; 
that is, that the sets $\m_j$ and $\m_{j'}$ do not intersect
nontrivially when $j\neq j'$.
We already know that this holds when $j$ and $j'$ are comparable.

Let $k$ be some vertex in $T$, and suppose that $i$ and $i'$ are two
distinct children of $k$.  
Among the neighbors of $k$, let $v$ be the one which lies one step
clockwise from $i$, and let $v'$ relate to $i'$ in the same manner.
From the properties of our tree, we know that $v$ and $v'$ are
ancestors of $k$.  Therefore they are comparable.  It is 
obvious that $v\neq v'$.  Assume without loss of generality that $v>v'$.
From the Inductive Lemma we know that $\m_k$ touches the left side
of $\m_v$, and that $\m_{i'}$ touches the left side of $\m_{v'}$.
Furthermore, in the proof of the lemma we have seen that $\m_i$
is contained in the region ${R_v^k}$ (In the notation of the lemma
$v$ is $v_0$ and $k$ is $v_e=i-1$).
Since $v'<v$, the set $\m_{v'}$ is not in that region, and therefore
neither is $\m_{i'}$, which touches it.  If $j\geq i$ and $j'\geq i'$,
then $\m_j$ is in ${R_v^k}$, while $\m_{j'}$ is not.
We deduce that $\m_j$ cannot intersect $\m_{j'}$.

If we now start with any two noncomparable vertices $j,\ j'$,
we let $k$ be their meet, let $i$ be the child of $k$ which is
$\leq j$, and similarly for $i'$ and $j'$. The argument of the
previous paragraph shows that $\m_j$ and $\m_{j'}$ do not intersect,
as required.  This shows that we have a packing.

It now follows that our packing has the required nerve, for the argument
above shows that it contains no edges which are not in $T$,
and the Inductive Lemma~7.2 shows that all the edges of $T$ are present in
the nerve.  The fact that the packing is nondegenerated follows
from Claim~6.4 and $y\in \br L$.  
This completes the proof of the Monster Packing Theorem.
\qed\enddemo

\remark{Remark 7.3} 
The Inductive Lemma actually implies that the packing produced,
in addition to having $T$ as its nerve, also agrees with
the orientation of $T$.  For example, the sets $\m_j$ which
touch $\m_b$ touch the left side of $\m_b$.
\endremark

\newsec\heading
	8. Packing Convex Bodies
\endheading

In this section we will apply the Monster Packing Theorem to prove
the existence of packings of convex bodies prescribed up to homothety.
We will start with a situation in which a few simplifying assumptions
have been made, and then move on to rid ourselves of the
nonessential assumptions.

\proclaim{8.1 Proposition}
Let there be given
a triangulation of the sphere,
$T$, a clockwise triangular
2-cell $\tri abc$, in $T$, 
and for each
vertex $v$ in $V$, the set of vertices of $T$, let there
be a prototype $P_v$, which is a set in the plane.
Assume that $P_a$ is the unbounded region determined by a smooth
simple closed curve, and each $P_v,\ v\neq a,$ is a
smooth, strictly convex compact set in the plane.
Assume that every pair of the three sets $P_a,\ P_b,\ P_c$ touch each other. 
Then there is a packing in the plane $Q=(Q_v:v\in V)$ whose nerve is
$T$ such that $Q_a=P_a,\ Q_b=P_b,\ Q_c=P_c$ and each
$Q_v$ is (positively) homothetic to $P_v$, for
$v\in V-\{a,b,c\}$. (See Fig.~\newfig).
\endproclaim
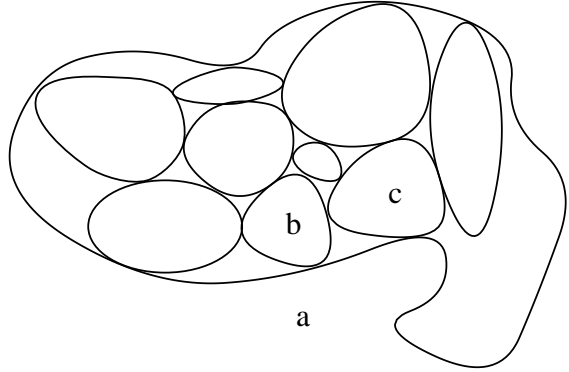
\fffffig \fig:A packing of convex shapes in a Jordan domain.:\fig:2in: 

Here, and in the following, a curve is called smooth,
if it has finite length and is once differentiable
with respect to an arc-length parameter.  A planar set is smooth, if
it has nonempty interior and its boundary is smooth.
 
\demo{Proof}
We will use the Monster Packing Theorem.  For that
purpose we will need to create a monster whose components
have the desired shapes.

Let $\leq$ be the order induced on $V$ by a right oriented
dfs tree of $T$ rooted at $\dedg ab$.
We will use the notation and terminology
introduced in the previous sections.
A monster $\m$ will be constructed in the following manner.
For every $x\in K=[0,1]^J,\ (J=V-\{a,b,c\})$,
set $\m_a(x) = P_a\cup\{\infty\}$,
$\m_b(x) = P_b$, $\m_c(x)=P_c$.  The remaining sets $\m_v(x)$ will be defined
inductively, using our partial order $\leq$.
Choose $p_b$ and $p_c$ to be some touching point of
$\m_a,\m_b$ and $\m_b,\m_c$ respectively and set
$p_b(x)=p_b,\ p_c(x)=p_c$.

Let $v\in J,\ x\in K$.
Assume that $\m_j(x),\ (j<v)$, and $p_j(x),\ (a<j<v),$ 
have been defined already, 
in such a way that the interiors of $\m_k(x),\m_j(x)$
are disjoint when $k<j$, that
every $p_j(x)$ is a point where $\m_j(x)$ touches $\m_{j-1}(x)$,
and that each $\m_j(x),\ a<j<v$, is either a point or is homothetic
to $P_j$.
If $\m_{v-1}(x)$ is a point, $\m_{v-1}(x)=\{p_{v-1}(x)\}$, and set 
$p_v(x)=p_{v-1}(x),\ \m_v(x)=\{p_v(x)\}$.
Otherwise let $p_v(x)$ be the point on the boundary of
$\m_{v-1}(x)$ so that the length of the part of the boundary
of $\m_{v-1}(x)$ extending counterclockwise from $p_{v-1}(x)$ to
$p_v(x)$ is exactly 
$x_v\cdot\operatorname{length}(\partial\m_{v-1}(x))$.
In other words, as $x_v$ increases from $0$ to $1$, while the other
coordinates remain fixed, the point $p_v(x)$ travels counterclockwise
along $\partial\m_{v-1}$, with constant speed, starting
and ending at $p_{v-1}$.  This defines $p_v(x)$.

Now there is not much choice left in defining $\m_v(x)$.
$\m_v(x)$ will be defined to be a set $S$,
which is either homothetic to $P_v$ or a point, touches
$\m_{v-1}(x)$ at $p_v(x)$, and touches $\cup _{j<v-1}\m_j(x)$.
(See Fig~\newfig).
We shall now see that such a set exists and is unique.
Since we assumed that
$\m_{v-1}(x)$ is not a point, it is homothetic to
$P_{v-1}$.  Since $\m_{v-1}(x)$ is convex and smooth,
a necessary and sufficient condition for a convex set $C$ to
touch $\m_{v-1}(x)$ at $p_v(x)$ is that the (unique) line $l$
which supports $\m_{v-1}(x)$ at $p_v(x)$ also supports $C$
at $p_v(x)$ on the other side of it.  
The set $P_v$ is strictly convex, therefore there is a unique point
$p$ on the boundary of $P_v$ such that the line parallel to $l$
passing through $p$ supports $P_v$ on the side opposite to
the side in which $l$ supports $\m_{v-1}(x)$.  Thus any set homothetic to
$P_v$ which touches $\m_{v-1}(x)$ at $p_v(x)$, is obtained from
$P_v$ by translating it by the translation which takes $p$
to $p_v(x)$ and then applying a homothety with center $p_v(x)$.
Conversely, any such set touches $\m_{v-1}(x)$ at $p_v(x)$ and
is homothetic to $P_v$.
Let $C$ be such a set.
First we will assume that none of the sets $\m_j(x),\ j<v-1$, contains
the point $p_v(x)$.
Let $C_\alpha, C_\beta$ be homothetic images of $C$ under homotheties
with center $p_v(x)$ and coefficients $\alpha>\beta>0$.
Since $C$ is strictly convex, $C_\alpha$ will contain the boundary of
$C_\beta$ in its interior, except for the point $p_v(x)$.
This shows that there is at most one 
positive number $\alpha$ such that the homothety with center 
$p_v(x)$ and coefficient $\alpha$ takes $C$ to a set
which touches $\cup_{j<v-1}\m_j(x)$.
Such a positive number exists, of course, because $\Bbb C-\m_a$ is bounded.
This establishes the existence and uniqueness in the case
$p_v(x)\notin \cup_{j<v-1}\m_j(x)$.
\fffffig \fig:Not much choice left in defining $\m_v(x)$.:\fig:4.2in: 

In case
$p_v(x)\in \cup_{j<v-1}\m_j(x)$, $p_v(x)$ is a point where
$\m_{v-1}(x)$ touches one of the sets $\m_j(x),\ j<v-1$.
None of these sets is a point, because otherwise
$\m_{v-1}(x)$ would also be, by the inductive construction.
It is impossible for three smooth sets to touch at a point,
and therefore the only possibility in this case is that
$\m_v(x)=\{p_v(x)\}$.  This verifies the existence and uniqueness
claimed above and completes the construction.

Now we need to verify that the mapping constructed is a monster,
that is, it satisfies conditions \ri{M1}--\ri{M7}.
Conditions \ri{M1}, \ri{M3}, \ri{M5} 
and the first assertion of \ri{M4} are immediate.
\ri{M7} follows from the fact that our sets are smooth.

Continuity will be proved by induction.
Let $v\in J$ and
assume that the sets $\m_j,\ j<v,$ and the points
$p_j,\ a<j<v,$ depend continuously on $x\in K$.
It is obvious that $p_v(x)$ is also continuous.
Let $x^1,x^2,\dots$ be a sequence in $K$ converging to the
point $x$.  By the inductive assumption $\m_j(x^i)$ converges
to $\m_j(x)$, as $i\to\infty$, for $j<v$.  We can choose a subsequence
of $(x^i)$, which we shall denote
by $(y^i)$, such that the sets $\m_v(y^i)$ converge to some set $S$.
For continuity it suffices to prove that $S=\m_v(x)$.

From the corresponding facts about the sets $\m_j(y^i),\ j\leq v,$
we know that, in the limit as $y^i\to x$, the set $S$ touches
one of the sets $\m_j(x),\ j<v-1$, its interior is disjoint
from these sets, and $S$ is either homothetic to $P_v$ or is a point.
Similarly, since $p_v(y^i)\to p_v(x)$ as ${i\to\infty}$,
$S$ touches $\m_{v-1}(x)$ at $p_v(x)$. 
$S$ shares these properties with $\m_v(x)$. They were the
properties defining $\m_v(x)$ in the case where
$\m_{v-1}(x)$ is not a point.  Thus $S=\m_v(x)$ in this case.

In the case that $\m_{v-1}(x)$ is a point, we have defined $\m_v(x)$
to be the same point, $\m_v(x)=\{p_v(x)\}$, and must verify that $S$ is
also.  Let $j<v$ be the smallest vertex such that 
$\m_j(x)$ is a point.  Necessarily $\m_j(x)=\{p_v(x)\}$.
Thus, by the definition of $\m_j(x)$,
there is a set $\m_k(x),\ k<j-1,$ touching $\m_{j-1}(x)$
at $p_v(x)$.  Since three smooth sets cannot touch at a point,
we must have $S=\{p_v(x)\}=\m_v(x)$, establishing continuity.
This yields \ri{M2} and the remaining part of \ri{M4}.
\ri{M6} is also clear.

Now the Monster Packing Theorem can be applied to yield our proposition.
\qed\enddemo
 
Now that we proved the proposition, let's try to discard some of
the unnecessary hypotheses, leading to Thm.~8.3.

\subheading{Freeing $P_b$ and $P_c$}
An unpleasant feature of Proposition~8.1 is the nonsymmetry in
the hypotheses required for the sets $P_a,P_b,P_c$,
while the conclusion is completely symmetric in them---The
configuration $(P_a,P_b,P_c)$ can be extended to a
packing with the required properties.  Our current objective
is to arrive at a more symmetric form for Proposition~8.1,
by dropping unnecessary requirements from $P_b,P_c$.
This more symmetric form will be convenient for applications.

The first fact to be observed is that the convexity of $P_b$ wasn't
used anywhere in the proof of the proposition; $P_b$ may be any smooth
set which is topologically a disk (and touches $P_a$ and $P_c$,
of course).  The convexity of $P_c$, on the other hand, was strongly
used for the construction of the monster in the proof of the proposition,
and therefore it will take some effort to dispense with it.
How was the convexity of $P_c$ used?
When a vertex $d$ of $T$ is a child of $c$, the point $p_d(x),\ x\in K$,
was defined to travel around the boundary of $P_c$ and the set
$\m_d(x)$, which is always homothetic to $P_d$ or a point,
touched $P_c$ at $p_d(x)$ and touched $P_a\cup P_b$.
It is easy to find examples which show that this isn't always possible
when $P_c$ isn't convex. 
To avoid this problem
we will need to modify the construction of the monster.

The construction given in the proof of the proposition for
the sets $\m_v(x)$, can remain exactly the same when $v-1>c$.
The only problem arises for the vertices $d$ in $T$ such that
$d-1=c$, that is, children of $c$.  
Since $\tri abc$ is a triangle in $T$, statement \ri{3}
of Prop.~4.1 implies that there is at most one child of $c$ in $T$;
denote it by $d$. (It is possible that $a,b,c$ are the only
vertices in $T$.  We ignore this trivial case).
So our task is to place $\m_d(x),\ x\in K$.
We examine the simpler case in which $P_d$ is a geometric disk, a
ball for the Euclidean metric of the plane.
The existence of a packing for the general case will follow from it.
 
Consider the following situation.  We are given three smooth sets
$P_a,\ P_b,\ P_c\subset\hat{\Bbb{C}}$ so that any pair of them touches.
Assume that the above sets are homeomorphic to disks and one
of them contains $\infty$ as an interior point.
We are interested in geometric disks, that is, closed balls of the
Euclidean metric on $\Bbb C$, which touch $P_c$ and touch
$P_a\cup P_b$.  Such disks we will call {\it permissible
disks}, and we do not exclude from this definition the
permissible disks of radius $0$, these are the sets containing
only a single point where $P_c$ intersects $P_a\cup P_b$.
A point on $\partial P_c$ which is contained in some
permissible disk will be called an {\it accessible point}.
The set of accessible points will be denoted by $A$.

\proclaim{8.2 Lemma}
$A$ is a closed subset of $\partial P_c$.
For each accessible point there is a unique permissible disk touching it.
If $\Gamma$ is an arc of $\partial P_c$ having both endpoints in
$A$, but no other point of it is in $A$, then the permissible disks
touching the two endpoints of $\Gamma$ are the same.
\endproclaim

\demo{Proof}
Let $q\in A$.  If $q$ is a point where $P_c$ touches $P_a$ or
$P_b$, then the only permissible disk containing $q$ is $\{q\}$,
because three smooth sets cannot touch at a point.
Assume that this is not the case.
Since $P_c$ is smooth,
any permissible disk containing $q$
must be tangent to $\partial P_c$ at $q$.  
Comparing two different disks which are tangent from the outside
to $P_c$ at $q$, we see that the larger one contains the boundary
of the smaller one in its interior, except for the point $q$.
Thus at most one of them can be a permissible disk, and the second
assertion of the lemma holds.

For the first assertion, let $(z_i)$ be a sequence of points in $A$
converging to some point $z\in\partial P_c$.
The corresponding sequence of permissible disks containing the
points $z_i$ will have some subsequence which will converge to
some disk $D$.  It is immediate that $D$ is a permissible
disk containing $z$.  Thus $z\in A$, as needed.

\comment
Now let $D_1,\ D_2$ be two distinct permissible disks with centers
$o_1,\ o_2$ and radii $r_1,\ r_2$ respectively.
Let $q_1$ be a point where $D_1$ touches
$P_a\cup P_b\cup P_c$, and similarly, let $q_2$ be such a point for $D_2$.
We will show now that the line segments $[o_1,q_1]$
and $[o_2,q_2]$ do not intersect.
Let $d_1,\ d_2$ be the distance of $q_2$ from $o_1$ and the distance
of $q_1$ from $o_2$, respectively.
Since $D_1\neq D_2$, the uniqueness argument above can be applied to show that
$d_1>r_1$ and $d_2>r_2$.
Assume that the segments $[o_1,q_1]$ and $[o_2,q_2]$ intersect at a point $p$.
Then the sum of the lengths of the segments $[o_1,p]$ and $[p,q_2]$
is at least $d_1$, and similarly,
$\text{length}[o_2,p]+\text{ length}[p,q_1]\geq d_2$.
But
$$
\multline
\text{length}[o_1,p]+\text{length}[o_2,p]
+\text{length}[q_1,p]+\text{length}[q_2,p]\\
=\text{length}[o_1,q_1]+\text{length}[o_2,q_2]
=r_1+r_2,
\endmultline
$$
and we get $r_1+r_2\geq d_1+d_2$; a contradiction.
This shows that the line segments $[o_1,q_1]$
and $[o_2,q_2]$ do not intersect.
\endcomment

Now let $D_1,\ D_2$ be two distinct permissible disks with centers
$o_1,\ o_2$ respectively.  Let $q_1$ be a point where $D_1$ touches
$P_a\cup P_b\cup P_c$ and similarly let $q_2$ be such a point for $D_2$.
The uniqueness argument above can be applied to show that $q_1\neq q_2$,
and we obviously have $o_1\neq o_2$.
The point $o_1$ is in the half-plane consisting of all the points
that their distance to $q_1$ is no larger than their distance
to $q_2$, and the respective statement holds for $o_2$.
This shows that the line segments $[o_1,q_1]$
and $[o_2,q_2]$ do not intersect.

Assume that the last assertion of the lemma fails.
Let the endpoints of $\Gamma$ be $q_1$ and $q_2$, let 
$D_1$ and $D_2$ be permissible disks touching at $q_1,\ q_2$,
let their centers be $o_1,\ o_2$, and let $s_1,\ s_2$ be points where
$D_1$ and $D_2$ touch $P_a\cup P_b$, respectively.
From the previous paragraph we know that the two curves 
$C_1=[q_1,o_1]\cup [o_1,s_1]$ and
$C_2=[q_2,o_2]\cup [o_2,s_2]$ do not intersect.
Take any simple curve $C_3$ contained in
$P_a\cup P_b$ and joining the points $s_1$ and $s_2$.
(The sets $P_a$ and $P_b$ intersect).
Look at the region disjoint from $P_c$ that is
bounded by $\Gamma,\ C_1,\ C_3,\ C_2$. (See Fig.~\newfig).
In this region there is some point $o_3$ whose distance to
$P_c$ is the same as its distance to $P_a\cup P_b$.  $o_3$
must be the center of a permissible disk $D_3$.  Let $q_3$ be
a point where $D_3$ touches $P_c$.  From the previous paragraph
we know that the line segment $[o_3,q_3]$ does not intersect
$C_1\cup C_2$.  Any interior point of it
cannot intersect $C_3\cup \Gamma$.
We may conclude that $q_3\in \Gamma$, in contradiction to
our assumptions.  This contradiction shows that $D_1=D_2$, and
completes the proof of the lemma.
\qed \enddemo

\fffffig \fig:Looking for another permissible disk.:\fig:4.5in: 

We will use the lemma to construct a monster for our new
situation.  Set $\m_v(x)=P_v$ for $v=a,b,c,\ x\in K$.
($\m_c$ will need to be slightly modified later.)
Take some parametrization of $\partial P_c$,
$q:[0,1]\to\partial P_c$, which circles counterclockwise
around $P_c$ and such that $q(0)=q(1)=p_b$, and $q$ is
one-to-one otherwise.  
With each $t\in [0,1]$ we will associate a permissible
disk $D(t)$, in the following manner.  If $q(t)$ is
an accessible point, let $D(t)$ be the unique permissible disk
touching $q(t)$.  Otherwise, since $A$ is a closed set,
there is a unique arc $\Gamma$, which satisfies the hypotheses
of the lemma and contains $q(t)$.  In this case, we let $D(t)$
be the permissible disk which touches the endpoints of $\Gamma$.
Finally, let $\m_d(x)=D(x_d)$, for $x\in K$.  
The lemma shows that the definition makes sense.  

Is $\m_d(x)$ continuous in $x$?
If $q(t)$ is not an accessible point, then $D(t')=D(t)$ for
every $t'$ close enough to $t$, so $D(t)$ is continuous at such
$t$.  The uniqueness of the permissible disk touching
at a given accessible point shows that $D(t)$ is also
continuous at points $t$ such that $q(t)$ is accessible.
Thus $D(t)$ and hence also $\m_d(x)$ are continuous.

For this to become part of our monster, we must define
a point $p_d(x)$ which varies continuously, and is always
a point where $\m_d(x)$ touches $\m_c(x)$.  Shortly we will
see that this is not always possible with our current definitions
of $\m_c,\ \m_d$, and therefore we will slightly modify $\m_c$.
When $q(x_d)$ is an accessible point leave $\m_c(x)=P_c$, as
before, and set $p_d(x)=q(x_d)$.
Let $0\leq t_1<t_2\leq 1$, and assume that $q(t_1),q(t_2)$ are accessible,
but no point $q(t),\ t_1<t<t_2$, is accessible.
Let $\Gamma$ be the image of $[t_1,t_2]$ under $q(\cdot)$,
then $\Gamma$ is as in the lemma.
As $x_d$ increases from $t_1$ to $t_2$,
the set $\m_d(x)$ remains the permissible disk, say $D$, which
touches $\m_c$ at $q(t_1)$ and $q(t_2)$, but $p_d(x)$ must somehow
travel continuously from $q(t_1)$ to $q(t_2)$, while continually
remaining in $\m_c\cap\m_d(x)$.
To perform this difficult task we redefine $\m_c(x)$.

One of the arcs of the boundary of $D$ which extends between
$q(t_1)$ and $q(t_2)$ bounds, together with $\Gamma$, a region
$R$, which is disjoint from $P_a\cup P_b\cup P_c\cup  D$.  Denote
this arc by $C$.  
Take some smooth homotopy $H$ from $\Gamma$ to $C$ which
stays in the region $R$ and keeps the endpoints $q(t_1),q(t_2)$
fixed.
Now define $\m_c(x)$, for $t_1\leq x_d\leq t_2$, as
follows.  When $x_d$ increases in the first third of the interval
$[t_1,t_2]$ the set $\m_c(x)$ is $P_c$ together with part of $R$ as
varies according to the homotopy $H$, until, when $x_d$ reaches 
the point $(2t_1+t_3)/3$,
the set $\m_c(x)$ becomes $P_c$ together with the
closure of $R$.
(To be more explicit, if the homotopy is such that $H_0=\Gamma$
and $H_1=C$, then let $\m_c(x)$ be the region bounded
by $(\partial P_c-\Gamma)\cup H_s$, when $x_d=t_1+s(t_2-t_1)/3$.)
When $x_d$ increases from 
$(2t_1+t_3)/3$ to $(t_1+2t_3)/3$, $\m_c(x)$ remains the same as for
$x_d=(2t_1+t_2)/3$.  In the final third of $[t_1,t_2]$, $\m_d(x)$
reverses its expansion of the first third, and returns, when
$x_d=t_2$, to being $P_c$.
Now $p_d(x)$ can accomplish the jump from $q(t_1)$ to $q(t_2)$.
As $x_d$ varies in the first third of $[t_1,t_2]$, let $p_d(x)$
remain at $q(t_1)$.  In the second third let it travels along $C$
from $q(t_1)$ to $q(t_2)$.  In the last third let it remain stationary
at $q(t_2)$.

It is clear that with these definitions $\m_c(x),\m_d(x)$ and $p_d(x)$
are continuous.  The definition of $\m_v(x),p_v(x),\ v>d$ can
remain as in the proof of the proposition. 
By our construction, the mapping $\m$ thus constructed satisfies
the requirements \ri{M1}--\ri{M7}, and is therefore a monster.
The Monster Packing Theorem may be applied to yield a packing
$Q=(Q_i:i\in V)$.
Is $Q$ the packing we need?  Almost.  There are only two problems.
The first one is that it may happen that $Q_c\neq P_c$.  This is
because we modified the definition of $\m_c$.  If this happens, and we
replace $Q_c$ by $P_c$, we will still have a packing,
because the difference between $Q_c$ and $P_c$ can only
be in the region $R$ of the previous paragraph, and, in this region
$R$, there aren't any bodies $Q_v$, since the nerve of $Q$ is $T$.
The new packing in which $Q_c$ is replaced by $P_c$ will, of course,
still have nerve $T$.

The second problem is that we have restricted ourselves to
the case that $P_d$ is a geometric disk.  To overcome this restriction
we insert a new vertex in the triangulation $T$.  In $T$ the child
of $c$, $d$, is characterized by the fact that $\tri cda$ is a
clockwise triangular 2-cell.  Insert in this 2-cell of $T$
a new vertex $d'$, connect edges from it to $a,c,d$, and split
the 2-cell accordingly.  Let the resulting triangulation be $T'$.
By the restricted case, there is a packing with nerve $T'$
for the data that we have together with the
requirement that the set corresponding to $d'$ is a disk.
If this set is ignored one has the required packing for $T$.

Let us summarize what we have done.  We have shown that in
Proposition~8.1 one can relax the requirements imposed on $P_a, P_b,$ and
$P_c$.
It is sufficient to require of them that they are smooth topological
disks in $\hat {\Bbb C}$, with each one touching the other two,
and that one of them contains $\infty$ in its interior.

Suppose that $P_a,P_b,P_c$ satisfy these conditions,
then there are precisely two connected components of
$\hat{\Bbb C}-(P_a\cup P_b \cup P_c)$ which neighbor with all three
of $P_a,P_b,P_c$, and for just one of these connected components, say $W$,
$P_a,P_b$ and $P_c$ are counterclockwise situated, in that order.
So the packing given by Prop.~8.1 will necessarily take place in $W$,
because $\tri abc$ is a clockwise triangle in $T$.
(See Remark~7.3).
Therefore only the parts of the
boundaries of $P_a,P_b,P_c$ which are exposed
to $W$ are relevant to the Proposition.
So suppose that $C$ is a simple closed curve in $\Bbb C$, which is smooth,
except for three distinct points, and that $P_a,P_b,P_c$ are the
three arcs determined by these points,  situated counterclockwise
with respect to $W$, the bounded region determined by $C$,
and further suppose that the three angles of $W$ that $P_a,P_b$ and $P_c$
determine at their intersections 
are zero, then the existence of a packing like that in
the Proposition holds in this case also.

\subheading{Replacing Strict Convexity With Convexity}
Using approximations, the requirement that the sets $P_v,\ v\neq a,b,c,$
are strictly convex is easy to replace by the weaker condition
that they are convex.

To see this, suppose that the sets $P_v,\ v\in V-\{a,b,c\},$ 
are smooth compact convex sets, without the assumption that they are
strictly convex.
Since compact convex sets can be approximated
by smooth strictly convex sets (The intersection of all balls of radius
$r$ containing $S$ is a good strictly convex
approximation for $S$, if $r$ is large
and $S$ is compact and convex.  The smoothness is no problem.),
for each $v\in V-\{a,b,c\}$ we can find a sequences $P_v^1,P_v^2, \dots$ of 
smooth strictly convex sets converging to $P_v$. 
By the proposition,
for $i=1,2,\dots$, there
is a packing $Q^i=(Q_v^i:v\in V)$ with the right nerve,
with $Q_a^i=P_a,\ Q_b^i=P_b,\ Q_c^i=P_c$ 
and with each $Q^i_v,\ v\in V-\{a,b,c\},$ homothetic to
$P^i_v$ respectively.
We can choose a subsequence of packings 
such that all the sets will converge.
In our carelessness, we shall denote this subsequence too by $(Q^i)$.
So each sequence $(Q_v^i)$ converges to some set $Q_v$.
How can $Q=(Q_v:v\in V)$ fail to be a packing as required?
It cannot fail.
First it is clear that the interiors of the sets are disjoint
and that each of them touches all the sets that it is supposed
to touch.  Once we convince ourselves that the sets do not
degenerate to points it will be clear that each $Q_v,\ v\in V-\{a,b,c\},$
is homothetic to $P_v$, respectively, and that we have a packing
as required.
   
Certainly the sets $Q_a,Q_b,Q_c$ are not points, because they
are equal to $P_a,P_b,P_c$ respectively.
To reach a contradiction, assume that $Q_v$ consists of one point
$p$.  Let $V'$ be the set of vertices $u\in V$ such that $Q_u=\{p\}$,
and let $U$ be the set of vertices in $V-V'$ neighboring
with some vertex in $V'$.
For every $u\in U$,
the set $Q_u$ contains the point $p$
(and is therefore not a single point and thus is homothetic
to $P_u$).  Since at most two smooth sets can touch at $p$, 
we see that there can be at most two vertices in $U$.
This is clearly impossible, because $U$ must contain
at least one cycle 
(the edges with both endpoints in $U$ separate the union
of 2-cells with at least one vertex in $V'$ from the triangle
$\tri abc$),
and there are no multiple edges in $T$. 
This contradiction establishes our claim that the strict
convexity in Proposition~8.1 may be replaced by convexity.
For this argument, showing that the limit
packing does not degenerate, I am indebted to Richard Schwartz.
Compare also \cite{Th, Ch.~13}.

\subheading{Smoothness}
The discussion above exposes very well the extent to which
smoothness is an essential hypotheses in the Proposition.
There is no problem to approximate convex sets by smooth
convex sets, or to approximate a nonsmooth region $P_a$ by
smooth regions.  The place where the above reasoning breaks
down in the absence of smoothness is the proof of nondegeneracy.
In fact, it is easy to see that this cannot be overcome, at
least not in the generality in which the proposition is stated.
Fig.~\newfig\ illustrates an example which shows that some smoothness
requirement is needed.
\fffffig \fig:Impossible to pack nondegenerately.:\fig:5.3in: 

If one still doesn't like the smoothness restriction, one 
has two options.  One can try to find weaker conditions which
still suffice, or can drop the restriction altogether, and
see what can still be proved in this new situation.  We will pursue
both trails.
For the first endeavor, note that in the approximation argument
above all that was used was that for every cycle $C$ in $T$
the sets corresponding to vertices in $C$ cannot all touch at one
point.
This is satisfied, for instance, if, in the situation where
$P_a,P_b,P_c$ are arcs of a simple closed curve (this situation is
discussed at the end of the subsection about freeing $P_b$ and $P_c$),
the sets $P_v,\ v\in V-\{a,b,c\}$ are smooth, the arcs $P_a,P_b,P_c$ are
smooth, and each of the three angles of $W$ which they determine is
$<\pi$.

Without the presence of a smoothness condition,
or something like the no-three-touch rule, the approximation
procedure might yield a \lq degenerated packing\rq.  It 
will still retain the structure of the packing but some
degeneration may occur.  Each set will still be disjoint from the
interiors of the other sets, and will touch those sets which it
is supposed to touch by the decree of the triangulation $T$.  However
it may also touch some of the sets without the consent of $T$,
and some of the sets might degenerate to points (or line segments,
if some of the $P_v$ are line segments).
We will term such a configuration a {\it degenerated packing conforming
to} $T$.

To summarize:

\proclaim{8.3 Convex Packing Theorem}
Let there be given a triangulation of the sphere,
$T$, and a clockwise triangular 2-cell $\tri abc$, in $T$.
Let $C$ be a simple closed curve in the plane, and 
let $P_a,P_b,P_c$ be three arcs composing $C$ in counterclockwise order,
which are determined by three distinct points of $C$.
For each vertex $v$ in $V-\{a,b,c\}$ ($V$ is the set of vertices of $T$),
let there be a prototype $P_v$, which is a convex set in the plane
containing more than one point.
Then there is a (possibly degenerated) packing in the plane 
$Q=(Q_v:v\in V)$ conforming to $T$, such that $Q_a=P_a,\ Q_b=P_b,\ Q_c=P_c$ 
and each $Q_v\ (v\in V-\{a,b,c\})$ is either a point or 
(positively) homothetic to $P_v$.

In the following situation it can be guaranteed that the packing
above is nondegenerate, and its nerve is $T$.
The sets $P_v\ (v\in V-\{a,b,c\})$
are also smooth, the arcs $P_a,P_b,P_c$ are smooth, and the three angles
of the bounded region determined by $C$ (at the points where these arcs meet)
are $<\pi$.\qed
\endproclaim

\subheading{Shape Fields}
We will now generalize Theorem~8.3 still further.
Instead of determining the sets to be packed up to homothety,
they will be specified in a more general manner.  This will be
useful later for some applications.

\definition{8.4 Definition}
Let $D$ be a set in the plane.
A (convex) {\it shape field\/} on $D$ is a continuous function $\f$
which assigns to each point $p\in D$
a compact convex set $\f(p)$ containing more than a single point.
The {\it image\/} of $\f$ 
is then the collection of all
sets $P\subset D$ which are homothetic to $\f(p)$ for some point $p$
on the boundary of $P$.
$\f$ will be called {\it smooth}, if all the sets $\f(p)$ are smooth.
(It is not required that $\f$ is a smooth function in any sense.)
Similarly, when the sets $\f(p)$ are strictly convex, $\f$ is said
to be {\it strictly convex}.
\enddefinition

We will see some examples of shape fields in the next section.
The reason for introducing these definitions is the following
theorem.

\proclaim{8.5 Shape Field Packing Theorem}
Let there be given a triangulation of the sphere,
$T$, and a clockwise triangular 2-cell $\tri abc$, in $T$.
Let $C$ be a simple closed curve in the plane, and 
let $P_a,P_b,P_c$ be three arcs composing $C$ in counterclockwise order,
which are determined by three distinct points of $C$.
Let $D$ be the union of $C$ and the bounded domain determined by $C$.
For each vertex $v$ in $V-\{a,b,c\}$ ($V$ is the set of vertices of $T$),
let $\f_v$ be a shape field on $D$.
Then there is a (possibly degenerated) packing in the plane 
$Q=(Q_v:v\in V)$ conforming to $T$, such that $Q_a=P_a,\ Q_b=P_b,\ Q_c=P_c$ 
and each $Q_v\ (v\in V-\{a,b,c\})$ is in the image of $\f_v$.

If the shape fields $\f_v$ are smooth,
the arcs $P_a,P_b,P_c$ are smooth, and the three angles
of the bounded region determined by $C$ are $<\pi$,
then the packing above is nondegenerate, and its nerve is $T$.
\endproclaim

\demo{Proof}
First the shape fields $\f_v$ are extended continuously to the whole
plane.  There is no problem in doing so, since $D$ is closed.

Consider now the statement analogous to Prop.~8.1, with the
requirement that the packed sets $Q_v,\ v\in V,$
be in the corresponding images of smooth, strictly convex
shape fields $\f_v$ replacing
the requirement that they be homothetic to some prototype, $P_v$.
The proof of this modified statement is the same as the proof of
Prop.~8.1, the only minor difference being in the inductive construction
of the monster $\m$.  When $\m_{v-1}$ and a point $p_v$ on the
boundary of $\m_{v-1}$ have been chosen, the set $\m_v$ is then
chosen as the only set homothetic to $\f(p_v)$ which touches
$\m_{v-1}$ at $p_v$ and also touches $\cup_{j<v-1}\m_j$.\footnote
{This idea, using the base $p_v$ to determine the homothety type,
is due to Richard Schwartz.}
The rest of the proof proceeds like the proof of Prop.~8.1.

Theorem~8.5 follows from the modified version of Prop.~8.1
in the same way that Theorem~8.3 followed from the original version.
\qed\enddemo

\newsec\heading 
	9. Applications to Conformal and Quasiconformal Mappings\\
	Of Multiply Connected Domains
\endheading

\proclaim{9.1 Theorem}
Let $G$ be an $n+1$-connected bounded domain in $\Bbb C$
which is obtained from a simply connected region $H$
after $n$ disjoint compact
connected sets $F_1,F_2,\dots,F_{n}$ have been removed from it.
Assume that none of the $F_j$-s is a point (to avoid trivialities).
Then for every planar bounded simply connected region $H'$,
and for every list 
$P_1,P_2,\dots P_{n}$ of $n$ compact convex sets which are not points,
there are disjoint sets
$P_1',P_2',\dots,P'_{n}$ contained in $H'$, with
each $P_j'$ (positively) homothetic to $P_j$,
such that $G$ is conformally equivalent to 
$
H'-\cup_{j=1}^{n}P_j'
$
and the boundary of each $F_j$ corresponds to
the boundary of the respective $P_j'$ under such a conformal equivalence.

In the case where the boundaries of $H$ and $H'$ are simple
closed curves, $\gamma $ and $\gamma '$ respectively, if three distinct
counterclockwise ordered points are chosen $z_1,z_2,z_3\in \gamma $,
and similarly  $z_1',z_2',z_3'\in \gamma '$,
then such $P_j'$-s and such an equivalence can be chosen
to satisfy the additional requirement that
each $z_j$ corresponds to $z'_j$ when the conformal
equivalence is extended continuously to $\gamma ,\gamma '$.
\endproclaim

The two ingredients entering the proof will be the
Convex Packing Theorem (8.3) 
and the method used in \cite{R-S} to prove
 that circle packings can yield
approximations of conformal mappings.
We use common results about quasiconformal mapping, all of which may be
found in \cite{L-V}, for instance.

Before the proof, let us introduce some notation.
In this section, when we speak of a circle, we shall mean the circle
together with its interior; that is, a closed disk.
A {\it hexagonal circle packing\/} is an infinite packing
of circles all having the same diameter so that each circle touches
precisely six other circles.  
It is unique, up to similarity.
The flower of a circle in such a packing consists of the circle
itself, its six immediate neighbors, and the six triangular
regions which they determine. 

\demo{Proof}
We deal with the situation in the second part of the theorem,
and make the following additional assumptions.  We assume that
each $P_j$ is smooth and strictly convex. Denoting by
$\gamma'_1$ the arc of $\gamma'$ extending counterclockwise from $z_1'$ to $z_2'$,
by
$\gamma'_2$ the arc of $\gamma'$ extending counterclockwise from $z_2'$ to $z_3'$,
and by
$\gamma'_3$ the arc of $\gamma'$ extending counterclockwise from $z_3'$ to $z_1'$,
we further assume that $\gp_1,\gp_2,\gp_3$ are smooth and that the angles
of $H'$ they determine at $z_1',z_2',z_3'$ are $<\pi$.
It is not hard to see, and will also become clear from our proof,
that once we have verified the theorem under these additional hypotheses,
the extra hypotheses can be removed through approximation.
Thus it is sufficient to restrict ourselves to this situation.

Let $\gam_1,\gam_2,$ and $\gam_3$ be the arcs of $\gam$ determined
by $z_1,z_2,$ and $z_3$,
as $\gp_1,\gp_2,$ and $\gp_3$ where determined by $z_1',z_2',$ and $z_3'$.
Choose some point $p$ in $G$.
For every small $\eps>0$, let $H_\eps$ be a hexagonal circle
packing of the whole plane with the diameters of the circles being $\eps$.
Provided $\eps$ is small enough, the flower of one of these circles, say $C_0$,
will contain $p$, and be contained in $G$.
A circle of $H_\eps$ will be called {\it inner\/} if it can be connected
to $C_0$ by a chain of circles whose flowers are contained in $G$.
A circle of $H_\eps$ will be called a {\it boundary\/} circle
if it is not inner, but is adjacent to some circle which is.
Let $Q=Q_\eps$ be the packing consisting of the inner and boundary
circles in $H_\eps$, and let $N=N_\eps$ be its nerve.
See Fig.~\newfig\  
and Fig.~\newfig. 
We take $\eps$ small enough so that the collection of boundary
circles will consist of $n+1$ disjoint cycles of circles, $R,S_1,\dots,S_n$:
the circles in $R$ will be those which are near $\gamma$, and those
in each $S_j$ are near $F_j$, respectively.
{
  \fffffig \fig:The packing $Q_\eps$ approximates the region $G$.:9.2:4in:
}
\fffffig \fig:The nerve $N$.:9.3:4.5in: 

Split $R$ into three paths of circles, $R_1,R_2,R_3$, with an overlap
of one circle between every two of these paths, so that
the distance between any circle of $R_j$ to $\gam_j$ is at most
$2\eps$, for $j=1,2,3$.  It is clear that this can be done.

We shall construct from $N$, the nerve of the packing $Q_\eps$,
a triangulation $T=T_\eps$.  
In (the complement of the embeding in the sphere of)
$N$ there are 2-cells $R^*,S_1^*,\dots,S_n^*$ whose vertices
correspond to the circles in $R,S_1,\dots,S_n$ respectively.
For each $j=1,2,\dots,n$ take a new vertex, $v_j$, insert it
into the 2-cell $S_j^*$ and split $S_j^*$ by inserting edges
which join $v_j$ to each vertex in the boundary of $S_j^*$.
Into $R^*$ insert three vertices, $a_1,a_2,a_3$, connect them with edges
to form a clockwise triangular 2-cell, $\tri {a_1}{a_2}{a_3}$,
connect $a_1$ to all the vertices corresponding to circles in $R_1$,
connect $a_2$ to all the vertices corresponding to circles in $R_2$,
connect $a_3$ to all the vertices corresponding to circles in $R_3$,
and split $R^*$ into triangular 2-cells accordingly.
One easily verifies that this procedure yields a triangulation of the
sphere.  See Fig.~\newfig.
\fffffig \fig:The triangulation $T$.:9.4:8in: 

By the Convex Packing Theorem, there exists some packing $H_\eps'$ whose
nerve is $T$, and which satisfies the following.  
The sets in $H_\eps'$ corresponding to
$a_1,a_2,a_3$ are $\gp_1,\gp_2,\gp_3$, respectively, those corresponding to 
$v_1,v_2,\dots,v_n$ are homothetic to $P_1,P_2,\dots,P_n$, respectively,
and the rest, those corresponding to vertices which are in $N$, 
are circles.
We shall denote the set in $H_\eps'$ corresponding to $v_j$ by $P'_{j,\eps}$;
it is homothetic to $P_j$.  
$Q'_\eps$ will denote the collection of circles in $H_\eps'$
corresponding to the vertices which were originally in $N$.
So there is a one-to-one correspondence between the circles of $Q_\eps$
to those of $Q_\eps'$.

Pick some sequence of $\eps$-s decreasing to zero.
The remainder of the proof uses
the methods of \cite{R-S} to show that some subsequence
of the \lq approximate mappings\rq\, which match each circle in the packing 
$Q_\eps$ to the circle corresponding to the same vertex of $T$ in $Q_\eps'$, 
converges to a conformal mapping, as required.

\proclaim{9.2 Resolution Lemma}
The diameters of the circles in $Q'_\eps$ decrease to $0$ as $\eps\to 0$.
\endproclaim
\demo{Proof of Resolution Lemma}
The proof is an adaptation of the proof of the Length-Area lemma of \cite{R-S}.
Let $d$ be a positive number which is smaller than 
the distance between any two of the sets $\gamma,F_1,F_2,\dots,F_n$.
Let $C$ be some circle in $Q_\eps$, let $C'$ be the corresponding
circle in $Q_\eps'$, and
let $k$ be the largest integer satisfying $(k+2)\eps<d/2$.
Look at the $k$ concentric rings of circles in $H_\eps$ around $C$ 
consisting of $6,12,18,\dots,6k$ circles, denote these rings
by $Z_1,Z_2,\dots,Z_k$.  ($Z_j$ consists
of the circles in $H_\eps$ that have combinatorial distance $j$ from $C$).
They meet at most one of the sets $\gamma,F_1,F_2,\dots,F_n$.
Consider first the case in which they don't meet any of these sets.
In this case they are in $Q_\eps$.  Denote the corresponding
rings of circles in $Q_\eps'$ around $C'$ by $Z_1',Z_2',\dots,Z_k'$.
Let $l_j$ denote the length of $Z_j'$, that is, 
twice the sum of the radii of circles in $Z_j'$, and let
$l$ be the minimum of the $l_j\ (j=1,2,\dots,k)$.
In $Z_j'$ there are $6j$ circles and the sum of their radii is $l_j/2$.
From the Cauchy-Schwarz inequality, we deduce that the sum
of the squares of the radii is $\geq (l_j/2)^2/6j$.  This means that the area
covered by all the $Z_j'$-s together is
$\geq\text{constant}\cdot l\sum_{j=1}^k 1/j$,
because the $Z_j'$ are disjoint in area.
But this area is certainly bounded above by the area of $H'$, thus,
as $\eps\to 0$ and $k\to\infty$, necessarily also $l\to 0$.
But if $l$ is small then $C'$ is also small, because each $l_j$ is
an upper bound for the diameter of $C'$.  This is so, because the closed
curve consisting of line segments between consecutive circles 
in $Z_j'$ has length $l_j$ and surrounds $C'$.

We deal now with the case where the concentric rings $Z_1,Z_2,\dots,Z_k$ 
meet $\gamma$
(and therefore do not meet any of the sets $F_1,F_2,\dots,F_n$).
Since these rings are not all contained in $Q_\eps$, we need to modify the
argument of the previous case.
We can assume that $d$ has been chosen small enough so that there is no
point in the plane whose distance to all three arcs $\gam_1,\gam_2,\gam_3$
is $\leq d$.
Assuming this, one of these arcs, say $\gam_1$, has distance $>d$ from
the center of $C$.
Thus, each circle in one of the rings $Z_1,Z_2,\dots,Z_k$ has
distance $>2\eps$ from $\gam_1$, and is not in $R_1$.
Denoting now by $Z_j'$ the collection of circles in $Q_\eps'$ which
correspond to circles of $Q_\eps$ which are in $Z_j$, we see that each
$Z_j'$ is either a ring of circles separating $C'$ from 
$\gp_1$, or contains a path
of circles which, together with portions of $\gp_2,\gp_3$, separates
$C'$ from $\gp_1$.
Define $l$ as twice the minimum of the sum of radii of circles
in $Z_j'$, as previously. 
So in $H_\eps'$ there is either a
closed curve of length $l$ separating $C'$ from $\gp_1$,
or a path of length $\leq l$, which together with portions
of $\gp_2$ and $\gp_3$ does so.
If $l$ is small, this implies that $C'$ is small.
(The following fact was used here.
If $x,y$ are points on the simple curve $\gam_2\cup\gam_3$, 
then the diameter of the portion of the curve between $x$ and
$y$ tends uniformly to zero,
as the distance between $x$ and $y$ tends to zero).
Applying the same proof as in the previous case, we know that 
$l$ tends to zero when $\eps\to 0$.
(If $Z_j'$ contains less circles than $Z_j$,
it makes our estimates even better).
So we conclude that in this case also $C'$ is small once $\eps$ is.

The case where the rings $Z_1,Z_2,\dots,Z_k$ intersect one of the sets
$F_1,F_2,\dots,F_n$ is dealt with similarly, but less care is needed.  
The only difference is that the $Z_j'$-s and portions of one of the
$P'_{m,\eps}$ now separate $C'$ from $\gam$; that is, they surround $C'$.
This concludes the proof of the lemma.
\qed \enddemo

We now wish to define a mapping, by utilizing the correspondence between
$Q_\eps$ and $Q'_\eps$.  Later, it will be shown that this mapping
is \lq approximately conformal\rq\ when $\eps$ is small.  
For a packing whose sets are all circles, the nerve has a natural
geometric realization.  In this realization, the vertices are the centers of
the appropriate circles, and an edge between vertices is the line segment
joining them.  The 2-cells are, of course, the regions determined by
these edges, polygons.
So we have two geometric realizations for the nerve $N$, one coming from
the packing $Q_\eps$, and the other, from $Q_\eps'$.  
This is the basis of the definition of our \lq approximately conformal\rq\ 
map, $f_\eps$.  

Look at the 2-cells of $N$. Except for $R^*,S_1^*,S_2^*,\dots,S_n^*$,
every such 2-cell is a triangular 2-cell, and has a geometric realization 
coming from $Q_\eps$, as a geometric triangle $A$,
and a geometric realization, coming from $Q_\eps'$, as a geometric
triangle $A'$.
Such a triangle $A$ will be called an {\it inner triangle\/} (of $Q_\eps$), if
its vertices correspond to inner circles of $Q_\eps$.
$G_\eps$ will denote the union of such inner triangles.
The triangle $A'$, corresponding to an inner triangle $A$ of $Q_\eps$,
will be called an {\it inner triangle\/} of $Q'_\eps$.
The union of the inner triangles of $Q'_\eps$ will be denoted by $G_\eps'$.

Let $f_\eps$ be the mapping which is defined on $G_\eps$,
takes each vertex of an inner triangle $A$ to the corresponding
vertex of the corresponding triangle $A'$, and is affine in each such $A$.
In other words, $f_\eps$ takes the center of an inner circle of $Q_\eps$ to
the center of the corresponding circle of $Q_\eps'$ and is
affine in each of the triangles which have as vertices the centers of three
inner circles of $Q_\eps$ that touch.  
It is immediate that $f_\eps$ is well-defined, and is a homeomorphism
between its domain, $G_\eps$, and its range, $G_\eps'$.

The following assertions are evident:
Every $G_\eps$ is contained in $G$,
and every compact subset of $G$ is contained in $G_\eps$,
once $\eps$ is small enough.
$G_\eps$ is the union of equilateral triangles, of edge-length $\eps$.
The boundary of $G_\eps$ is contained in a $3\eps$ neighborhood of
the boundary of $G$.

\proclaim{9.3 Quasiconformality Lemma}
There is some constant $K$ such that each mapping $f_\eps$ is
$K$-quasiconformal. 
Furthermore, given a compact
subset $B\subset G$ and a number $\alpha>1$, then $f_\eps$ is 
$\alpha$-quasiconformal on $B$, once $\eps$ is small enough.
\endproclaim
\demo{Proof}
If seven circles form a flower, as in Fig.~\newfig,
then it is easy to see
that the ratio of the radius of any one of the surrounding circles
to the center circle is bounded away from zero.  (This can be
deduced from the fact that three circles cannot touch
at a point.  Also see \cite{R-S,Ring Lemma}). 
This means that the ratio of the radii of two touching circles
of $Q_\eps'$ which correspond to inner circles of $Q_\eps$ is bounded.
We may conclude that the ratio of any two edges of an inner triangle
of $G_\eps'$ is bounded.
Take an inner triangle $A$.  $A$ is equilateral and $f_\eps$ maps it affinely
onto a triangle which is not too far away from being equilateral.
This implies that there is some constant $K$ so that every $f_\eps$ is
$K$-quasiconformal on every inner triangle,
and from this the first assertion of the lemma follows.
\fffffig \fig:Seven circles form a flower.:9.5:2.5in: 

In \cite{R-S} and also in \cite{He1}, it has been proven that the hexagonal
packing is rigid in the following sense.
In any circle packing which is combinatorially equivalent to $k$
generations of the hexagonal packing around a circle $C$, the
ratio of the radius of $C$ to the radius of 
any of its neighbors approaches 1, if $k$ approaches $\infty$.
(\lq $k$ generations of the hexagonal packing around a circle $C$\rq\
is the collection of circles having combinatorial distance from $C$ which
is $\leq k$.)
Take a triangle $A$ of $G_\eps$ intersecting $B$.  
A circle of $Q_\eps$ which corresponds to a vertex in such a triangle
is surrounded in $Q_\eps$ by many generations of the hexagonal packing,
and $\text{many}\to\infty$ as $\eps\to 0$.
The same is true in $Q'_\eps$ for the corresponding circle.
Applying the above quoted result to the situation in $Q_\eps'$,
we deduce that $A'$, the triangle corresponding to $A$ in $Q'_\eps$,
is as close as we can wish to being equilateral, provided $\eps$ is
small enough.  So given any $\alpha>1$, once $\eps$ is small enough,
$f_\eps$ is $\alpha$-quasiconformal on any of the triangles composing 
$G_\eps$ which meet $B$.  Since $G_\eps$ contains $B$ for sufficiently
small $\eps$, this implies our assertion.
\qed \enddemo

The equi-quasiconformality just proved will help us deduce the following.

\proclaim{9.4 Separation Lemma}
Let $z$ be some point in $G$.
As $\eps\to 0$ the distance of $f_\eps(z)$ to 
$\cup_{j=1}^n P_{j,\eps}'\cup \gam'$ is bounded away from zero.
\endproclaim
\demo{Proof}
Let $\eta$ be some compact subset of $(\gp_1\cup\gp_2)-\{z_1', z_3'\}$.
We first show that the distance of $f_\eps(z)$ to $\eta$ is bounded
away from zero.  

We will construct a quadrilateral $B$ having edges $B_1,B_2,B_3,B_4$.
(Fig.~\newfig).
Let $z_4,z_5,z_6, z_7$ be four points placed counterclockwise on 
$\gam$ between $z_3$ and $z_1$; that is, in $\gam_3$.
Let $B_1$ be the subarc of $\gam_3$ between $z_4$ and $z_5$.
Let $B_2$ be a simple curve with endpoints $z_5,z_6$, which lies
in $G$, except for its endpoints, and separates $z$ from
$\gam_1\cup\gam_2$ in $G$.
Let $B_3$ be the subarc of $\gam_3$ between $z_6$ and $z_7$.
Finally, let $B_4$ be a simple curve with endpoints $z_7,z_4$, which lies
in $G$, except for its endpoints, so that the strip of $G$ between
$B_2$ and $B_4$ does not contain any of the sets $F_1,F_2,\dots,F_n$.
This strip is taken as the interior of our quadrilateral $B$.
It is clear that $B_4$ separates $B_2$ from $\gam_1\cup\gam_2$ in $G$.
\fffffig \fig:The quadrilateral $B$.:9.6:6in: 

We approximate the quadrilateral $B$ by a quadrilateral $B_\eps$ in $G_\eps$
which has sides $B_{\eps,1},\ B_{\eps,2},\ B_{\eps,3},\ B_{\eps,4}$.
We do this in such a way that the distance from a point on any side 
$B_{\eps,j}$ of
$B_\eps$ to the corresponding side $B_j$ of $B$ is at most $3\eps$,
and so that $B_{\eps,2}$ separates $z$ in $G_\eps$ from the points
of $G_\eps$ near $\gam_1\cup\gam_2$, and $B_{\eps,4}$ separates $B_{\eps,2}$
in $G_\eps$ from the points of $G_\eps$ near $\gam_1\cup\gam_2$.
(It is easy to see that this can be done; one can take $B_{\eps,1}$
and $B_{\eps,3}$ to lie on the boundary of $G_\eps$).
As $\eps\to 0$ the quadrilateral $B_\eps$ tends to $B$, 
and the modulus of $B_\eps$ approaches the modulus of $B$. 
Since the $f_\eps$ are equi-quasiconformal, by the above lemma,
we can deduce that the moduli of the quadrilaterals
$B_\eps' \overset\text{def}\to=f_\eps(B_\eps)$ are bounded away from
zero and infinity.  Set $B_{\eps,j}'=f(B_{\eps,j})$.

What do we know about the sides of the quadrilaterals $B_\eps'$?
Since $B_{\eps,1}$ 
did not stray more than a constant times $\eps$ away from $\gam_3$,
we know that any point on $B_{\eps,1}'$ is at most
a constant number of circles of $Q_\eps'$ away from $\gp_3$.
A similar argument holds for $B_{\eps,3}$.
From the Resolution Lemma, we may therefore deduce that as $\eps\to 0$,
the sides $B_{\eps,1}'$ and $B_{\eps,3}'$ are confined to smaller and
smaller neighborhoods of $\gp_3$.

Suppose that the distance of $\eta$ to $\gp_3$ is $d$.
Let $\eps$ be so small so that the distance from $\eta$ to
$B_{\eps,1}'\cup B_{\eps,3}'$ is larger than $d/2$.
Suppose that the distance of $f_\eps(z)$ to $\eta$ is $\delta$,
and that $\delta<d/4$.
Since the quadrilateral $B_\eps'$ separates $f_\eps(z)$ from the
points of $G_\eps'$ that are near $\eta$,
we see that in the quadrilateral $B_\eps'$ there is a curve
$c$ of length smaller than $\delta$ which connects the sides
$B_{\eps,2}$ and $B_{\eps,4}$, and whose distance to
$B_{\eps,1}'\cup B_{\eps,3}'$ is $>d/4$.
Since the modulus
of $B_\eps'$ is bounded away from $0$ and $\infty$,
this implys that $\delta$ cannot be arbitrarily small
in relation to $d$.  This establishes our claim that
the distance of $f_\eps(z)$ to $\eta$ is bounded
away from zero.

The same argument can yield stronger results.  The
quadrilateral $B$ could be chosen to separate an arbitrary compact
subset of $G$ from $\eta$.  Therefore,
the images under $f_\eps$ of such a set are bounded away from $\eta$.
The same is true, of course, if we choose $\eta$ to be a compact subset
of $(\gp_2\cup\gp_3)-\{z_2',z_1'\}$ or a compact subset of
$(\gp_3\cup\gp_1)-\{z_3',z_2'\}$.
Since $\gp$ can be represented as a union of a compact subset of
$(\gp_1\cup\gp_2)-\{z_1', z_3'\}$,
a compact subset of $(\gp_2\cup\gp_3)-\{z_2',z_1'\}$, and a compact
subset of $(\gp_3\cup\gp_1)-\{z_3',z_2'\}$, we see that the images under
$f_\eps$ of a compact subset of $G$ are bounded away from $\gp$,
as $\eps\to 0$.

Now let $\Gam_1$ be a simple closed curve in $G$ which separates $z$
from $\cup_{j=1}^nF_j$, let $\Gam_2$ be a simple closed curve passing
through $z$, which separates $\Gam_1$ from $\gam$.
The images of $\Gam_1$ under the $f_\eps$ are bounded away from $\gp$.
Therefore, by looking at the image under $f_\eps$ of the annuli bounded by
$\Gam_1$ and the boundary component of $G_\eps$ corresponding to $\gam$,
and using the equi-quasiconformality of the $f_\eps$ again,
we see that the diameter of $f_\eps(\Gam_1)$ is bounded away from
zero.
Now looking at the images of the annulus bounded by $\Gam_1$
and $\Gam_2$, we can therefore conclude that the distance
of $f_\eps(\Gam_1)$ to $f_\eps(\Gam_2)$ is bounded away from zero,
because the moduli of these annuli are bounded away from zero.
Since $f_\eps(z)\in f_\eps(\Gam_2)$, and $f_\eps(\Gam_1)$ separates
$f_\eps(z)$ from $\cup_{j=1}^nP'_{j,\eps}$, we see that $f_\eps(z)$ is
bounded away from $\cup_{j=1}^nP'_{j,\eps}$.
This completes the proof of our lemma.
\qed \enddemo

\demo{Conclusion of the Proof of Theorem 9.1}
Pick some sequence $\{\eps_k\}$ converging to zero, so that for each
$j=1,2,\dots,n$ the sets $P'_{j,\eps_k}$ converge,
as $k\to\infty$, to some set, which will be denoted by $P_j'$.  
Let $G'=H'-\cup_{j=1}^nP_j'$.  From the Resolution Lemma it follows that
any compact subset of $G'$ is contained in the image of $f_{\eps_k}$,
once $k$ is large enough. 
The equi-quasiconformality implies
that none of the sets $P_j'$ is a point, that
these sets are disjoint, and therefore that $G'$ is connected.
Because of the equi-quasiconformality,
by picking a subsequence if necessary, we may assume that the functions
$f_{\eps_k}$ converge uniformly on compact subsets of $G$, and
their inverses converge uniformly on compact subsets of $G'$.
Assume this, and denote the limit function of $f_{\eps_k}$ by $f$,
and the limit function of $f^{-1}_{\eps_k}$ denote by $g$.
Take some point $z\in G$.  The Separation Lemma implies that the
set $\{f_{\eps_k}(z):k>k_0\}$ is contained in a compact subset
of $G'$, once $k_0$ is large enough.  Therefore, since the functions
$f^{-1}_{\eps_k}$ converge uniformly to $g$ on compact subsets of $G'$,
we have
$$
g(f(z))=\lim_{k\to\infty}f^{-1}_{\eps_k}(f_{\eps_k}(z))=z,
$$
and $f$ is a homeomorphism.
From the second part of the Quasiconformality Lemma we see
that $f$ is 1-quasiconformal, and therefore conformal.

To see that $f$ maps $G$ surjectively onto $G'$, it is sufficient to
show that $g(G')\subset G$, because then we can write
$f(g(w))=\lim_{k\to\infty}f_{\eps_k}(f^{-1}_{\eps_k}(w))=w$.
Clearly, $g(G')$ is contained in the
closure of $G$, so it is sufficient to demonstrate that $g$ is an
open mapping.
For this we quote the following result (see \cite{L-V, p.~74}):
The limit of a locally uniformly converging sequence of
$K$-quasiconformal maps defined in an open connected set is either a
constant, or a $K$-quasiconformal map, and therefore a homeomorphism.
Let $w_1,w_2$ be distinct points in $f(G)$ and let $U$ be
any connected open set containing $w_1,w_2$ whose closure is
contained in $G'$.  The function $g$ is not constant in $U$,
it takes distinct values in $w_1,w_2$.
Therefore, by applying the above quoted result to the sequence
$f_{\eps_k}^{-1}$ restricted to $U$, we conclude that $g(U)\subset G$.
This clearly implies $g(G')\subset G$, and $f$ is onto $G'$.

Each $P'_j$ cannot be a point, because $f$ is conformal.
Therefore, since it is the limit of the $P'_{j,\eps_k}$-s,
it must be homothetic to $P_j$.

It is clear that the boundary of each $F_j$ corresponds under $f$ 
to the respective $P_j'$.
Since $f$ can be (uniquely) extended continuously to $\gam$,
it remains to verify that this extension maps $\gam_1,\gam_2,\gam_3$
to $\gp_1,\gp_2,\gp_3$ respectively.
To see this, recall the quadrilateral $B$
from the proof of the Separation Lemma.
There we have shown that the images
under the $f_\eps$-s of the region of $G$ 
which is separated by $B$ from $\gam_1\cup\gam_2$
stay bounded away from any compact subset of $\gp_1\cup\gp_2$.
Since $B$ can be chosen to separate an arbitrary point of
$\gam_3-\{z_2,z_3\}$
from $\gam_1\cup\gam_2$,
we see that $f$ maps points of $\gam_3-\{z_2,z_3\}$ to points of
$\gp_3$.
Similar statements are true for $\gam_1,\gp_1$ and $\gam_2,\gp_2$.
This establishes the boundary value requirement, and completes the proof.
\qed\enddemo
\enddemo

\subheading{A Generalized Beltrami Equation}
Theorem~9.1 gives the existence of a map satisfying the Cauchy-Riemann
equations along with some boundary conditions.  We will now present
an analogue of Thm.~9.1 with equations more general than the Cauchy-Riemann 
equations.
In the situation of the theorem, let $z$ vary in $\br G$,
and let $w$ vary in $\br{H'}$.  Suppose that to each such pair $(z,w)$ 
corresponds an ellipse $E(z,w)\subset{\Bbb C}$ in a continuous manner.
The ellipse  $E(z,w)$ is allowed to have arbitrary eccentricity and
orientation. 
(The position and size of $E(z,w)$ will be irrelevant for our purposes.)

\proclaim{Theorem 9.5} Let the situation be as in Theorem~9.1, 
and let $E(z,w)$ be an ellipse field as above.
Theorem~9.1 is still true if instead of requiring
the existence of a conformal mapping,
we require the existence of a quasiconformal
mapping $f:G\to H'-\cup_{j=1}^nP_j'$ with the property that for almost
every point $z\in G$ the inverse image of $E(z,f(z))$ under
the differential of $f$ is a circle:
$$
df_z^{-1}E(z,f(z))\text{ is a circle},\qquad\text{for almost every } z\in G.
\tag 9.1
$$
The requirements other than conformality stay the same.
\endproclaim
Theorem~9.1 is obtained as a special case when one takes
all the ellipses $E(z,w)$ to be circles.


If $g=g(w)$ denotes the inverse of $f$, then condition \thetag{9.1} can
be transformed to the form
$$
\frac{\partial g}{\partial \bar w}=\mu(g(w),w)\frac{\partial g}{\partial w}
\qquad\text{almost everywhere},
$$
where $\mu$ is a complex continuous function, $|\mu(z,w)|<1$, which encodes
the shape of the ellipse field $E(z,w)$.
Thus \thetag{9.1} can be thought of as a generalized Beltrami equation for
the inverse of $f$.

There is yet another way to express condition \thetag{9.1}.  Let
$h$ be a homeomorphism of $G$ onto $G'\subset H'$, then $h$ induces
the continuous ellipse field $E(h^{-1}(w),w)$ on $G'$.
This ellipse field induces a new conformal structure on $G'$, in the
following manner.  To specify a conformal structure one needs to
specify conformal coordinate charts.  Consider an open topological disk $D$
in $G'$.  Solutions to the Beltrami equation produce (an essentially unique)
quasiconformal homeomorphism from $D$ onto the open unit disk 
with the property that its differential at almost every point
$w$ takes the ellipse
$E(h^{-1}(w),w)$ to a circle.  Let this be a coordinate chart for
$D$.  Since a quasiconformal homeomorphism whose differential
at almost every point takes a circle to a circle is conformal,
it follows that these charts give a conformal structure for $G'$,
induced by $h$.

Condition \thetag{9.1} and the quasiconformality of $f$ are together
equivalent to the statement that $f$ is conformal, when we use for the
image of $f$ the new conformal structure which $f$ induces on it via
the ellipse field $E(f^{-1}(w),w),\ w\in f(G)$.

We shall now sketch the modifications needed in the proof of  
Theorem~9.1 to prove Theorem~9.5.
The construction of $Q_\eps$ and $T_\eps$ remains as above,
but the packing $H_\eps'$ will change.
Its combinatorics will remain the same, but
for every circle $C$ of $Q_\eps$ the corresponding set $C'$ in $H_\eps'$,
instead of being a circle,
will be an ellipse $C'$ which is homothetic to $E(z,w)$ for some
points $z\in C$ and $w\in C'$.
(Rather vaguely, one can say that the correspondence $C\to C'$ is
an approximate solution of our differential equation.)
To prove the existence of such a packing $H_\eps'$, we choose
an arbitrary point $z\in C$, for every circle $C$ of $Q_\eps$.
The functions $w\to E(z,w)$ are shape fields on $\br{H'}$,
and Theorem~8.5 can be used to guarantee the existence of the
required packing $H_\eps'$.

Now that we've got $H'_\eps$, the corresponding analogue of the
Resolution Lemma~9.2 is proved in exactly the same manner as the
proof given;
one only has to observe that there is a bound on the eccentricity
of the ellipses involved.  This is because the $E(z,w)$ vary
continuously as $(z,w)$ varies in a compact set.

The definition of the mapping $f_\eps$ has to undergo
some cosmetic modifications, because the contact point of
two touching ellipses might not lie on the line connecting
the centers of them.  For the new definition of $f_\eps$,
we use the following geometric triangulation induced by the
packing $H'_\eps$.  The vertices of this geometric triangulation
are the centers of the ellipses and the contact points of any two
touching ellipses.  The center point of each ellipse is joined
by a line segment to every contact point of it.  Two contact
points of an ellipse are joined by a line segment if, and only if,
they are consecutive contacts of the ellipse; i.e., there isn't
any other contact points on the shorter arc of the ellipse which
joins them.  It is easy to see that this procedure gives a
geometric triangulation (of part of $H'$).

The corresponding procedure is done for the packing $Q_\eps$, and one
gets a correspondence between triangles in both geometric triangulations.
Now the mapping $f_\eps$ is defined similarly to its definition in the 
original proof, but using these new, more refined, triangulations.

The first assertion of the Quasiconformality Lemma~9.3 and the
Separation Lemma~9.4 readily follow.
Here, again, all that is needed is the fact that the ellipses have
bounded eccentricity.  Then it is clear that some sequence of $\eps\to 0$
can be chosen so that the $f_\eps$ converge uniformly on compact
subsets.  If the limit function is denoted by $f$, then $f$
induces a new conformal structure on its image, the one
induced by the assignment
of the ellipse $E(f^{-1}(w),w)$ to each point $w\in f(G)$.
By continuity of the ellipses $E(z,w)$, it follows that near a
point $w$, for sufficiently small $\eps$ (and $f_\eps$ near $f$),
all the ellipses of $H'_\eps$ will be of homothety type close
to $E(f^{-1}(w),w)$.  This means that the second part
of the Quasiconformality Lemma holds, when one takes the new
conformal structure induced on the image of $f$, and we see that
$f$ is conformal, in this conformal structure; verifying the claim.

\subheading{Shape Fields, Revisited}
Another generalization of Theorem~9.1 is obtained if instead of
having the sets $P'_j$ be specified up to homothety, they are
required to be in the image of given shape fields.  For the
convenience of the reader we will state an abridged version of this
generalization.

\proclaim{9.6 Theorem} 
Let $H$ and $H'$ be bounded Jordan domains in $\Bbb C$, and
let their boundaries be the Jordan curves $\gamma$ and $\gamma'$,
respectively.
Set $G=H-\cup_{j=1}^n F_j$, where the $F_j$ are some disjoint
connected compact subsets of $H$ containing more than one point.
Let $\f_1,\dots,\f_n$ be shape fields on $\br{H'}$.
Suppose that three distinct counterclockwise ordered points
are chosen $z_1,z_2,z_3\in \gamma $,
and similarly  $z_1',z_2',z_3'\in \gamma '$.
Then there are sets $F_1',\dots,F_n'$, with each $F_j'$ in
the image of the shape field $\f_j$
(see the definitions of the previous section), and there is a conformal
equivalence of $G$ and $H'-\cup_{j=1}^n F_j'$, so that under
this conformal equivalence each $F_j$ corresponds to $F_j'$,
and each $z_k$ corresponds to $z_k',\ (k=1,2,3)$.
\endproclaim

\demo{Proof} The proof of Theorem~9.1 generalizes directly to prove
Theorem~9.6; the only change is that an appeal to Theorem~8.5 replaces
the corresponding appeal to Theorem~8.3.
\qed\enddemo

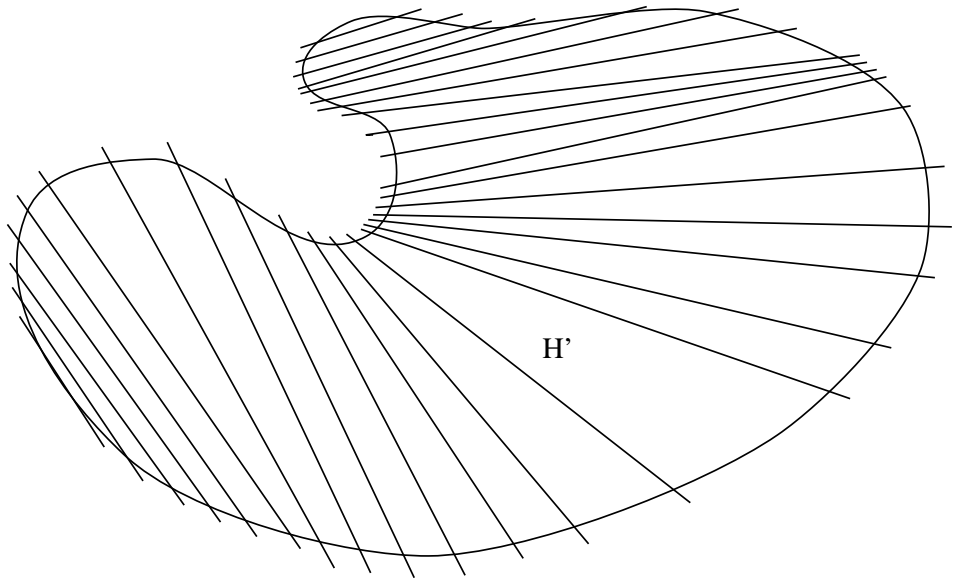
\fffffig \fig:A foliation of a neighborhood of $H'$ by line segments.:9.7:3.0in: 

We close this section with an example of a
class of shape fields having interesting images.
Let $H'$ be a Jordan domain.  Suppose that a neighborhood of $\br{H'}$ is
foliated by straight line segments (see Fig.~\newfig). 
One
can define a shape field on $\br{H'}$ by assigning to
each point $p\in \br{H'}$ a line segment of length one with
center at 0 whose direction is the direction of the foliation
at $p$.  The image of this shape field is the collection of
all line segments contained in the leaves of the foliation.
Given a finitely connected planar domain $G$,
one can use Theorem~9.6 and this example to construct
a conformally equivalent \lq slot domain\rq, with slots contained
in the leaves of the foliation.  This yields a generalization
of some of the canonical slot domains introduced by Koebe.
(See \cite{Si}).

\enddocument